\documentclass{article}

\usepackage[a4paper, total={6.05in, 8.2in}]{geometry}

\usepackage[latin1]{inputenc}   
\usepackage{amsmath}            
\usepackage{epstopdf}           
\usepackage{flushend}           
 \usepackage{setspace}

\usepackage{amsfonts, amsmath, calc,amssymb}
\usepackage{algorithm}
\usepackage[noend]{algpseudocode}

 \usepackage{setspace}
\let\Algorithm\algorithm
\renewcommand\algorithm[1][]{\Algorithm[#1]\setstretch{1.05}}

\usepackage{verbatim}
\usepackage{overpic}
\usepackage{comment}
\usepackage{tikz, tikz-cd}
\usetikzlibrary{arrows,shapes,cd}
\usepackage{caption}
\newcommand\bx{{\mathbf  x}}

\newcommand\bs{\sigma}

\newcommand\ba{\boldsymbol \alpha}
\newcommand\be{\boldsymbol e}

\newtheorem{note}{Note}
\newtheorem{remark}{Remark}

\newcommand*{\mc   }[2]{\multicolumn{#1}{c }{ {#2} } }
\DeclareMathOperator*{\argmin}{arg\,min} 

\definecolor{par_1}{rgb}{0.20784313725490197, 0.164705882352941,0.5294117647058824}
\definecolor{par_2}{rgb}{0.058823529411764705, 0.36078431372549,0.8666666666666667}
\definecolor{par_3}{rgb}{0.07058823529411765, 0.490196078431372,0.8470588235294118}
\definecolor{par_4}{rgb}{0.027450980392156862, 0.61176470588235,0.8117647058823529}
\definecolor{par_5}{rgb}{0.08235294117647059, 0.694117647058823,0.7058823529411765}
\definecolor{par_6}{rgb}{0.34901960784313724, 0.741176470588235,0.5490196078431373}
\definecolor{par_7}{rgb}{0.6470588235294118, 0.7450980392156863,0.4196078431372549}
\definecolor{par_8}{rgb}{0.8823529411764706, 0.7254901960784313,0.3215686274509804}
\definecolor{par_9}{rgb}{0.9882352941176471, 0.807843137254902, 0.1803921568627451}
\definecolor{par_37}{rgb}{0.3953, 0.7459, 0.5244}
\usepackage{amssymb}
\usepackage[figuresright]{rotating}
\usepackage{subfigure}
\usepackage{titling}
\usepackage{blindtext}
\usepackage[affil-it]{authblk}

\date{plain}
\begin{document}

\title{Efficient parallel optimization for approximating CAD curves featuring super-convergence}

\author[1]{Julia Docampo S\'anchez}

\affil[1]{Computer Applications in Science and Engineering,  Barcelona Supercomputing Center, Barcelona 08015, Spain}

 \maketitle
\begin{abstract}
We present an efficient, parallel, constrained optimization technique for approximating CAD curves with super-convergent rates.
 The optimization function is a disparity measure in terms of a piece-wise polynomial approximation and a curve re-parametrization. 
 The constrained problem solves the disparity functional fixing the mesh element interfaces.
  We have numerical evidence that the constrained disparity preserves the original super-convergence: ${2p}$ order for planar curves and $\lfloor\frac 32(p-1)\rfloor + 2$ for 3D curves, $p$ being the mesh polynomial degree. 
  Our optimization scheme consists of a globalized Newton method with a nonmonotone line search, and a log barrier function preventing element inversion in the curve re-parameterization. Moreover, we introduce a \emph{Julia} interface to the EGADS geometry kernel and a parallel optimization algorithm. We test the potential of our curve mesh generation tool on a computer cluster using several aircraft CAD models.  We conclude that the solver is well-suited for parallel computing, producing super-convergent approximations to CAD curves.
\end{abstract}

\section{Introduction and Motivation}

Geometric accuracy plays a major role in the performance of unstructured high-order methods \cite{cfd2030}, requiring curved elements to meet the desired accuracy. Traditionally, the geometric accuracy was measured based on the Fr\'echet and Hausdorff distances \cite{alt1995}. More recently, distance optimization techniques have been proposed. For example, in  \cite{remacle2014, toulorge2016} an area and Taylor-based distance optimizer are used respectively for 2D and 3D geometry. The authors report significant mesh-CAD distance reductions at adequate computational times.

In addition, a disparity measure for generating optimal curved high-order meshes was proposed in  \cite{ruiz2015}. The optimization combines a distortion measure for mesh quality and a geometric $L^2$-disparity measure for geometric error \cite{ruiz2016}. It produces optimal non-interpolative meshes and it has been observed that this disparity is $2p$ super-convergent \cite{ruiz2021}. This affords a straightforward advantage: one can obtain the desired geometric accuracy using smaller polynomial degrees than with standard interpolation approaches.

Solving the disparity implies solving a non-linear optimization problem. The excessive number of iterations limits the practical applications of this method. To overcome this, a constrained disparity with fixed element interfaces during optimization was proposed in \cite{docampo_imr}. Further, they employed a type of Zhang-Hager \cite{zhang} nonmonotonic line search based on the average history of the objective function. In general, nonmonotone line searches improve the computational efficiency as well as the likelihood of finding a global minimum \cite{zhang,grippo1996,toint, dai}. 

The work presented here focuses on practical applications of the disparity measure. It extends the constrained optimization approach from \cite{docampo_imr} to general CAD models and parallel processing. The optimization combines a globalized Newton method with the Zhang-Hager line search, and a log barrier function avoiding element inversion.
We have developed an interface to the EGADS geometry kernel \cite{egads} written in the \emph{Julia} language. The \emph{Julia} feature is now part of the ESP software distribution \cite{esp}, allowing direct use on clusters using the EGADSlite environment \cite{egadslite}  and \emph{Julia}'s HPC tools. 

Our results show that the constrained disparity reduces the dimension of the optimization problem, leading to fewer non-linear iterations. In addition, we run the optimization in parallel, further reducing the computational times. We test our methodology on several CAD models, running the simulations in parallel on one node of the \emph{BSC Marenostrum 4} Supercomputer.

This paper is organized as follows. In Section \ref{sec:disp}, we define the disparity measure for curves and show the super-convergent property through an example. In Section \ref{sec:new_solver}, we discuss the constrained disparity optimization including the algorithm parallelization. In Section \ref{sec:local_superconvergence} we discuss superconvergence for the constrained disparity in a global and local setting. Finally, in Section  \ref{sec:numerical_results}, we compare the errors and iterations from the constrained and unconstrained problem, and study the performance of the parallel scheme. The paper concludes with  Section \ref{sec:conclusion} discussing main results. 


\section{The original disparity measure}\label{sec:disp}

We begin by introducing the disparity formulation for curves, discussing how it is defined and optimized, as well as showing a superconvergent result. For more details, we refer the reader to \cite{ruiz2021}. 
\subsection{Mathematical formulation}
We define our mesh as a set of elements where for each physical element $e^P$ there is a reference element $e^R$. The physical mesh $\mathcal M^P$ can be defined in terms of an element-wise parametrization $\phi^P$:
\begin{align}\label{eq:x_mesh}
\phi^P|_{e^R}: e^R &\to e^P\subset \mathbb R^n\\
               \xi &\to \bx = \sum\limits_{i = 1}^{p+1} \bx_i N_i^p(\xi),
\end{align}
where $p+1$ are the number of nodes for the high-order element $e^P$, $\bx_i$ the $\mathbb R^n$-physical coordinates cooresponding to  the  $i$-th node and  $\{N_i^p\}_{i=1}^{p+1}$ a Lagrangian basis of degree $p$.

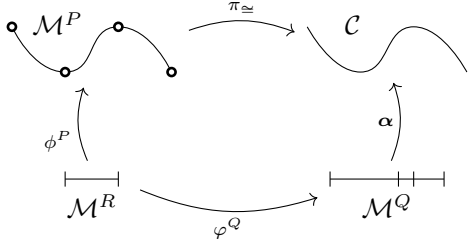
\begin{figure}[h!]
\centering
\begin{tikzcd}[ampersand replacement=\&,column sep=4em]
\begin{tikzpicture}
\draw (0.6,0.5) node {$\mathcal M^P$};
\draw  plot [ultra thick,smooth, tension=1.1] coordinates { (0,0.6) (0.7,0.0) (1.4,0.6) (2.1,0.0)};
\fill (0,0.6)        circle (2pt);
\fill (0,0.6) [white]circle (1pt);
\fill (0.7,0.0)        circle (2pt);
\fill (0.7,0.0) [white]circle (1pt);
\fill (1.4,0.6)        circle (2pt);
\fill (1.4,0.6) [white]circle (1pt);
\fill (2.1,0.0)        circle (2pt);
\fill (2.1,0.0) [white]circle (1pt);
\end{tikzpicture}
 \arrow[bend left=20]{r}{\pi_{\cong}} 
   \& 
\begin{tikzpicture}
\draw (0.6,0.5) node {$\mathcal C$};
\draw  plot [ultra thick,smooth, tension=1] coordinates { (0,0.6) (0.7,0.0) (1.4,0.6) (2.1,0.0)};
\end{tikzpicture}
\\[10pt]
\begin{tikzpicture}
\draw (0.7,   0)--(1.4,  0);
\draw (0.7,-0.1)--(0.7,0.1);
\draw (1.4,-0.1)--(1.4,0.1);
\draw (1.05,-0.5) node{$\mathcal M^R$};
\end{tikzpicture}
\arrow[bend left=20]{u}{\phi ^P}
\arrow[bend right=20,swap]{r}{\varphi^Q} 
\&  
\begin{tikzpicture}
\draw (0,0)--(1.5,0);
\draw (0  ,-0.1)--(0  ,0.1);
\draw (0.9,-0.1)--(0.9,0.1);
\draw (1.1,-0.1)--(1.1,0.1);
\draw (1.5,-0.1)--(1.5,0.1);
\draw (0.75,-0.5) node {$\mathcal M^Q$};
\end{tikzpicture}
\arrow[bend right=20]{u}{\ba} 
 \end{tikzcd}
 \caption{\label{fig:disp_diagram} commutative diagram showing the reference mesh $\mathcal M^R$, the mappings to the physical meshes: $\phi^P$ and $\phi^Q $ respectively, the curve $\mathcal C$, its parametrization $\ba$, and the projection $\pi$.}
 \end{figure}

Let $\mathcal C\subset \mathbb R^n$  be a curve parametrized by $\ba:\text [a,b]\subset \mathbb{R}\to \mathcal C$ and  consider the family of mappings:  
$\Pi:=\left\{\pi \in \mathcal H^1(\mathcal M^P,\mathcal C),\ \pi \text{ diffeomorphism}\right\}.
$
The diagram in Figure \ref{fig:disp_diagram} shows the reference mesh $\mathcal M^R$, the physical mesh $\mathcal M^P$ consisting of three elements, the target curve $\mathcal C$ and the diffeomorphism $\pi$.

The disparity measure between the high-order mesh and the curve is the minimal projection error
\begin{align}\label{eq:disp1}
d(\mathcal M^P, \mathcal C) &=  \inf\limits_{\pi\in\Pi} \left(\int\limits_{\mathcal M^P} | \bx - \pi(\bx) |^2 d\bx\right)^{1/2},
\end{align}
where $| \cdot |$ denotes the Euclidean norm of vectors.
Defining the functional 
\begin{align}\label{eq:func1}
E(\bx,\pi) &= \int\limits_{\mathcal M^P} | \bx - \pi(\bx) |^2 d\bx = \int\limits_{\mathcal M^R} | \phi^P(\xi) - \pi \circ \phi^P (\xi)|^2 |\dot{\phi^P}(\xi)| d\xi \\
 &= || \phi^P - \pi \circ \phi^P ||^2_\bs,
\end{align}
we establish the following relation:   
\begin{align}\label{eq:disp2}
d(\mathcal M^P, \mathcal C)^2 = \inf\limits_{\pi\in\Pi} E(\bx, \pi).
\end{align}
Here, the sub-index $\bs$ denotes that it is an integral with weight $|\dot\phi ^P|$.

Consider any possible curve  reparametrization: $\ba \circ s$. As in \eqref{eq:x_mesh}, we define a 1D mesh through the mapping:
\begin{align}\label{eq:s_mesh}
\varphi^Q|_{e^R}: e^R &\to e^Q\subset \mathbb R\\
               \xi &\to s = \sum\limits_{i = 1}^{q+1} s_i \cdot N_i^q(\xi).
\end{align}
\begin{note}

The mesh $\phi^P$ is in the physical space whereas $\varphi^Q$ is a mesh in the parametric space. Modifying $\varphi^Q$ results in different curve parametrizations.  
\end{note}

As shown in Figure \ref{fig:disp_diagram}, we have that $\pi\circ \phi^P = \ba \circ \varphi^Q$. Hence, we can reformulate the problem: 
\begin{align}{\label{eq:disp_s}}
  E(\bx, \pi) &= E(\bx, s) = || \bx - \ba \circ s ||^2_{\bs} = \int\limits_{\mathcal M^R} | \phi^P(\xi) - \ba \circ \varphi^Q (\xi) |^2 | \dot{\phi^P}(\xi) | d\xi.
\end{align}
Therefore, the mapping $s:\mathcal M^R \to \mathcal M^Q$ needs to be a diffeomorphism, too. In Section \ref{sec:new_solver}, we show how to enforce this numerically. 

The disparity between the mesh and the curve is found optimizing $E$ as a function of $s$.  If we optimize $E$ for both $\bx$ and $s$, we obtain the mesh with optimal geometric accuracy according to the disparity measure. We define the optimal approximation as:
\begin{align}{\label{eq:disp2}}
\bx^\star,s^\star &= \argmin\limits_{\bx,s} || \bx - \ba\circ s ||^2_{\boldsymbol \sigma} = \argmin\limits_{\bx,s}  \int\limits_{\mathcal M^R} | \bx(\xi) - \ba\circ s(\xi) |^2 |\dot{\bx}(\xi)| d\xi.
\end{align}

\subsection{An example of super-convergence}

\begin{figure}[h!]
\flushleft
\vspace{10pt}
\hspace{50pt}\begin{overpic}[
width=keepaspectratio,width=0.55\textwidth, abs,unit=1mm]
{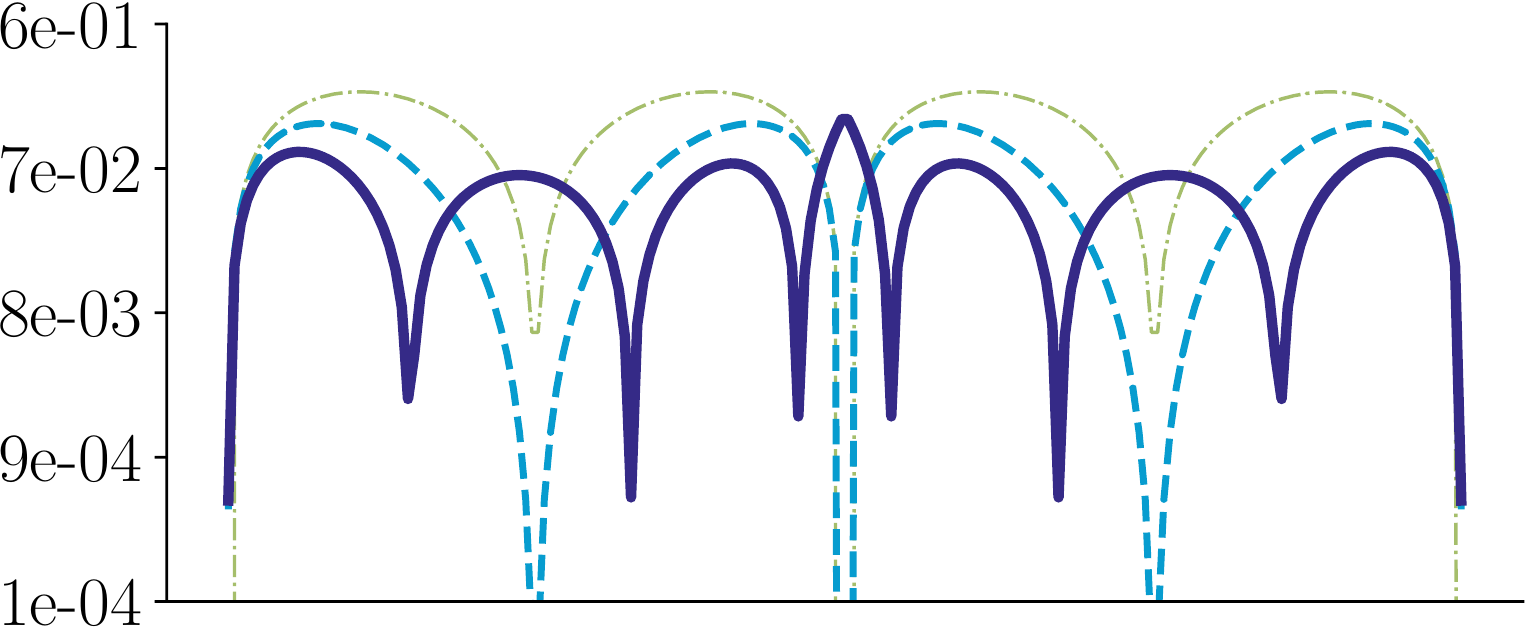}
\put(80,10){
\begin{tikzpicture}
\draw (0,1)--(0.4,1) [par_7, dashdotted   , line width = 1.0pt];
\draw (0,0.25)--(0.4,0.25) [par_4,densely dashed    , line width = 1.3pt];
\draw (0,-0.5)--(0.4,-0.5)     [par_1,                  , line width = 1.45pt];
\draw (1.2,1) node {$|\bx - \ba|$};
\draw (1.6,0.25) node {$|\bx - \ba\circ s^\star|$};
\draw (1.7,-0.5)   node {$|\bx^\star - \ba\circ s^\star|$};
\end{tikzpicture}}
\put(40,-2){\scriptsize{$\xi$}}
\put(11,-2){\scriptsize{$0$}}
\put(70,-2){\scriptsize{$2\pi$}}
\put(-10,10){\rotatebox{90}{$|error|$}}
\end{overpic}
\vspace{5pt}
\caption{\label{fig:circle_error} point-wise errors 
approximating a circle with two elements of degree $p=2$ obtained with: direct interpolation $|\bx(\xi)-\ba(\xi)|$ (dashed green), the disparity $|\bx(\xi) - \ba\circ s^\star(\xi)|$ (dotted light blue) and the optimized disparity $|\bx^\star(\xi)-\ba\circ s^\star(\xi)|$ (solid deep blue).}
\end{figure}
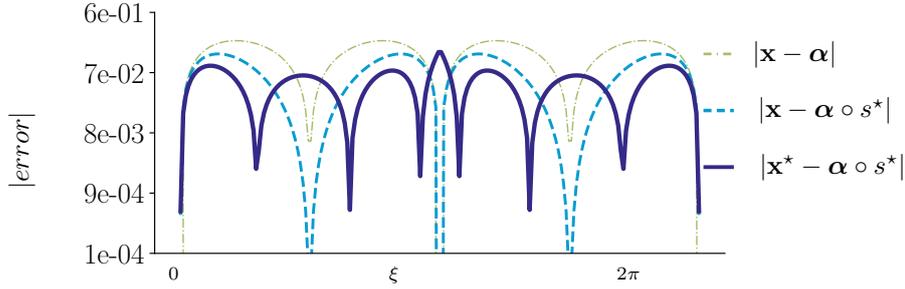

Let us discuss the role played by the physical  ($\bx$) and parametric ($s$) meshes. Figure \ref{fig:circle_error} shows several point-wise errors: first, we approximate the circle with interpolative meshes: $(\bx,s)$. Then, we compute the disparity measure of this mesh optimizing the parametric mesh,  $s^\star = \argmin\limits_{s} || \bx - \ba\circ s||^2_{\bs}$. Finally, we optimize both: $\bx^\star,  s^\star = \argmin\limits_{\bx,s} || \bx - \ba\circ s||^2_{\bs}$. 
Notice how the optimal mesh improves the geometric error of the original interpolation approximation.  

\begin{figure}[h!]
\flushleft
\hspace{40pt}
\begin{overpic}[
height=keepaspectratio,height=0.35\textwidth]
{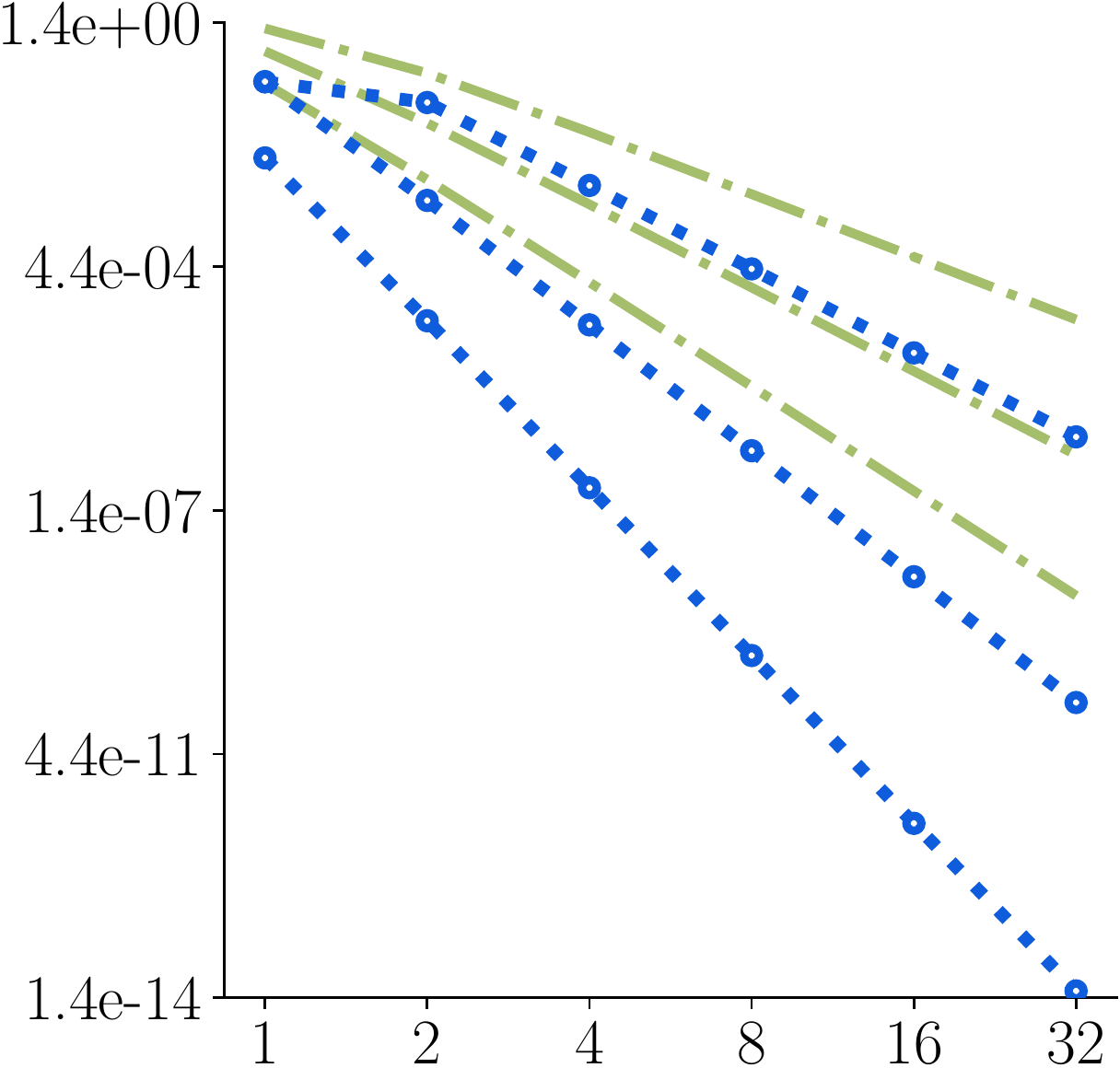}
\put(140,60){
\begin{tikzpicture}
\draw (0,1.0)--(0.5,1.0) [par_7, dashdotted, line width = 1.5pt];
\draw (0,0.0)--(0.5,0.0) [par_2,densely dashed, line width = 1.5pt];
\draw (1.6,1.0) node {$||\bx - \ba||_{\bs}$};
\draw (2.1,0.0) node {$||\bx^\star- \ba\circ s^\star||_{\bs}$};
\end{tikzpicture}
}
\put(66,33){
\tikz{\draw (0,0.2)--(0.4,-0.1)--(0.4,0.2)--cycle [par_2, line width = 1.pt]}}
\put(63,29){
\tikz{\draw (0,0.2) node{\textcolor{par_2}{$\scriptstyle{\boldsymbol{8}}$}}}
}

\put(66,55){
\tikz{\draw (0,0.2)--(0.4,-0.05)--(0.4,0.2)--cycle [par_2, line width = 1.pt]}}
\put(63,51){
\tikz{\draw (0,0.2) node{\textcolor{par_2}{$\scriptstyle{\boldsymbol{6}}$}}}
}

\put(66,75){
\tikz{\draw (0,0.2)--(0.4,0.05)--(0.4,0.2)--cycle [par_2, line width = 1.pt]}}
\put(63,70){
\tikz{\draw (0,0.2) node{\textcolor{par_2}{$\scriptstyle{\boldsymbol{4}}$}}}
}

\put(105,60){$p=2$}
\put(105,30){$p=3$}
\put(105,10){$p=4$}

\put(-15,45){\rotatebox{90}{$|| \cdot ||_{\bs}$}}

\end{overpic}

\caption{\label{fig:circle_rates_optimal} 
slopes (log-log) of the the initial norm $||\bx - \ba||_\sigma$ and optimized $||\bx^\star - \ba\circ s^\star||_\sigma$  disparities approximating a circle for several mesh refinements highlighting the super-convergent rates \protect\tikz{\protect\draw(0.0,0.3)--(0.4,0.1)--(0.4,0.3)--cycle[par_2, line width = 1.pt]} 
 .}
\end{figure}

In Figure \ref{fig:circle_rates_optimal}, we show convergence plots of the disparity measure when approximating  the same circle using meshes of degree $p=2,3,4$ and for several $h$-refinements. The initial disparity gives to $p+1$ order. On the other hand, the slope of the optimal pair ($\bx^\star,s^\star)$ shows a convergence order of $2p$. 

\section{Constrained disparity optimization and parallelization}\label{sec:new_solver}
Our new solver consists of four main ingredients: the \emph{constrained disparity} with fixed element interfaces during optimization (1), a \emph{Log barrier} function preventing curves from tangling (2), a globalized Newton method with the \emph{Zhang-Hager line search} (3), and \emph{parallel optimization} (4).

\subsection{The constrained disparity}
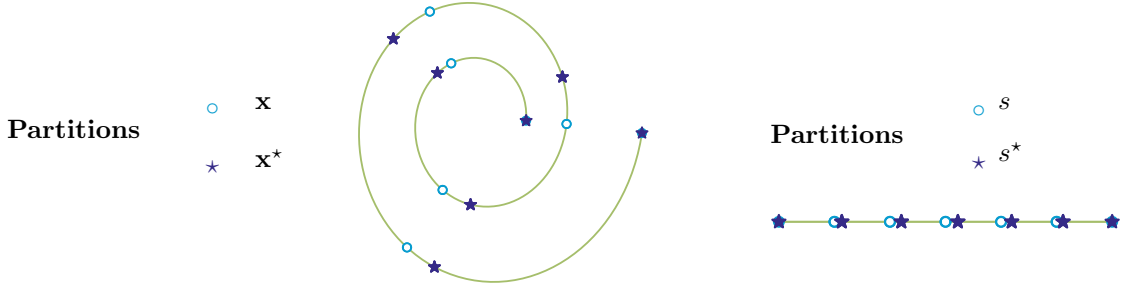
\begin{figure}[h!]
\flushright
\begin{tabular}{cc}
\begin{overpic}[width= keepaspectratio, width = 0.25\textwidth]
{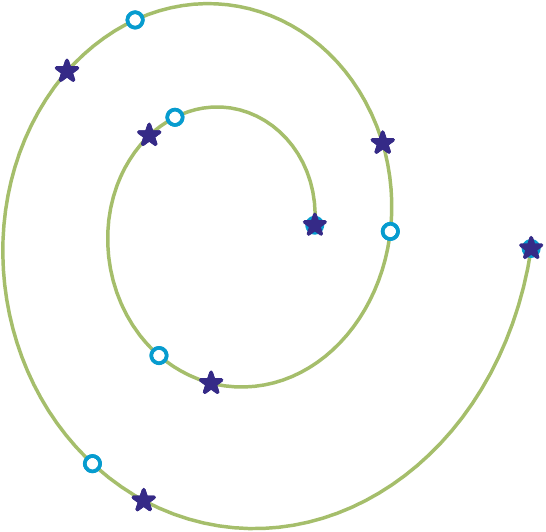}
\put(-120,50){\textbf{Partitions}}
\put(-35,60){$\bx$}
\put(-35,40){$\bx^\star$}
\put(-55,55){\tikz{\draw(0,0) node{\textcolor{par_4}{$\circ$}}}}
\put(-55,35){\tikz{\draw(0,0) node{\textcolor{par_1}{$\star$}}}}
\end{overpic}
&\hspace{30pt}
\raisebox{20pt}{\begin{overpic}[width= keepaspectratio, width = 0.3\textwidth]
{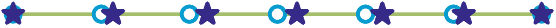}
\put(0,25){\textbf{Partitions}}
\put(65,35){$s$}
\put(65,20){$s^\star$}
\put(55,30){\tikz{\draw(0,0) node{\textcolor{par_4}{$\circ$}}}}
\put(55,15){\tikz{\draw(0,0) node{\textcolor{par_1}{$\star$}}}}
\end{overpic}}

\end{tabular}
\caption{\label{fig:xs_spiral_elements} a spiral approximated with 6 elements showing the physical ($\bx$) and parametric ($s$) element parititions (\textcolor{par_4}{$\circ$}, \textcolor{par_1}{$\star$} respectively) before and after optimizing the disparity.}
\end{figure}

In the previous section we defined our physical and parametric elements using a general map from a reference element (equations \eqref{eq:x_mesh} and \eqref{eq:s_mesh}). Since we consider all possible partitions along the curve, the element distribution changes as we optimize the disparity. In Figure \ref{fig:xs_spiral_elements} we illustrate this for a spiral curve; the initial and final element configuration is different in $\bx$ and $s$, with elements moving towards the spiral end.

We now cast a different disparity optimization problem by fixing the element interfaces. We give the mathematical formulation and highlight the differences concerning the original disparity. Recall that we optimize the geometric accuracy of a high-order mesh solving the problem:
\[
E(\bx^\star,s^\star) = \min\limits_{\bx,s} || \bx -\ba \circ s||^2_{\bs},
\]
with $\bx$ consisting of elements of degree $p$ and $s$ elements of degree $q$. 

Assume a fixed element partition and define each mesh element by
\begin{align}\label{eq:const_meshes}
 \tilde\bx(\xi) = \underbrace{\bx_0N_0^p(\xi) + \bx_pN_{p+1}^p(\xi)}_{\bx_F(\xi)} + \sum\limits_{i = 2}^{p}  \tilde \bx_i \cdot N_i^p(\xi), \\
 \tilde s(\xi)  = \underbrace{s_0N_0^q(\xi) + s_q N_{q+1}^q(\xi)}_{s_F(\xi)} +\sum\limits_{i = 2}^{q} \tilde s_i \cdot N_i^q(\xi),
\end{align}
where $\bx_0,\ \bx_{p+1},\ s_0,\ s_{q+1}$ are fixed throughout the optimization. Note that the indices run from 2 to $p$ and 2 to $q$ respectively, excluding the first and last nodes.
We define the optimal \emph{constrained} mesh as 
\begin{equation}\label{eq:disp_sub}
\tilde \bx^\star,\tilde s^\star = \argmin_{\tilde \bx,\tilde s} ||\tilde \bx - \ba\circ \tilde s||^2_{\bs}.
\end{equation}

\begin{note}
Uniform parametric partitions on the CAD may lead to poor initial element distributions, affecting the overall accuracy of the optimized meshes. 
When a curvature-based mesher is not available, we propose a pre-processing stage: find the optimal (in the disparity sense) element distribution optimizing the original (unconstrained) disparity functional with linear meshes ($\bx^\star_1,\ s^\star_1)$. The $s^\star_1$ will be the initial parametric partition on the CAD. Alternatively, one can perform arc-length-based optimization \cite{shontz2014}. 
\end{note}

\begin{figure}[h!]
\centering
\begin{tabular}{ccc}

\mc{3}{
\begin{tikzpicture}
\draw (-1.0, 0.0) node{\textbf{Curves:}};
\draw (0.5, 0.0)--(0.3, 0.0) [par_6,  densely dashed, line width =1pt];
\draw (1.2,0.0)  node{$\ba\circ \tilde s^\star$};
\draw (2.2, 0.0)--(2.5, 0.0) [par_1, line width =1pt];
\draw (3.1,0.0)  node{$\tilde \bx ^\star$};
\end{tikzpicture}
}\\[5pt]
\begin{overpic}
[width = keepaspectratio, width = 0.2\textwidth]{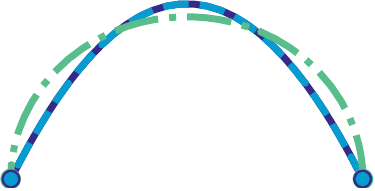}
\put(25,1){
\fbox{\textbf{R=1}}}
\end{overpic}
&
\begin{overpic}
[width = keepaspectratio, width = 0.2\textwidth]{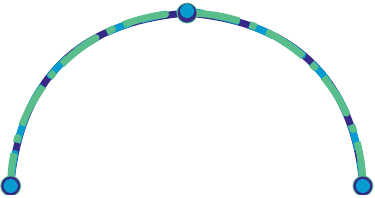}
\put(25,1){
\fbox{\textbf{R=2}}}
\end{overpic}
&
\begin{overpic}
[width = keepaspectratio, width = 0.2\textwidth]{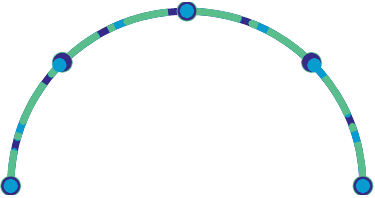}
\put(25,1){
\fbox{\textbf{R=4}}}
\end{overpic}
\\[15pt]
\mc{3}{
\begin{tikzpicture}
\draw (-2.2,0.0)  node{\textbf{Point-wise errors:}};
\draw (0.2, 0.0)--(0.6, 0.0) [par_6,  densely dashed, line width =1pt];
\draw (1.8,0.0)  node{$|\bx - \ba \circ s|$};
\draw (3.2, 0.0)--(3.6, 0.0) [par_1, line width =1pt];
\draw (5.0,0.0)  node{$|\tilde \bx ^\star - \ba \circ \tilde s^\star|$};
\end{tikzpicture}
}\\[5pt]
\includegraphics[
height=keepaspectratio,height=0.15\textwidth]{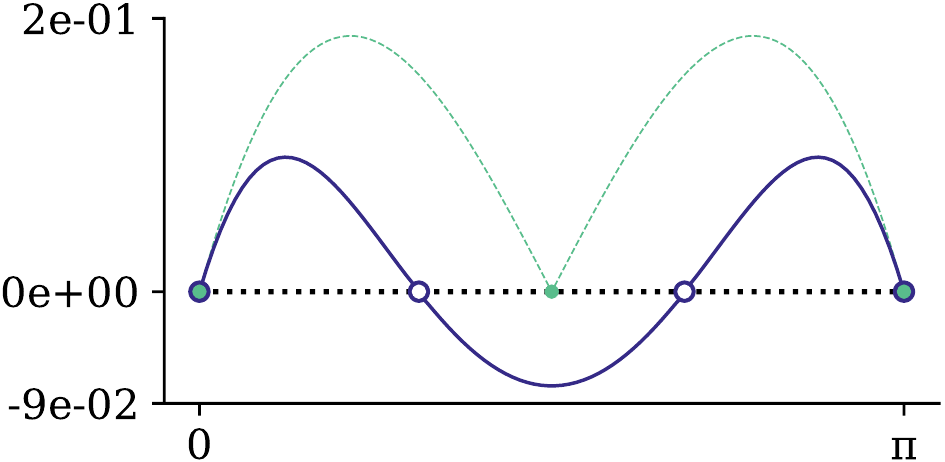} 
&
\includegraphics[
height=keepaspectratio,height=0.15\textwidth]{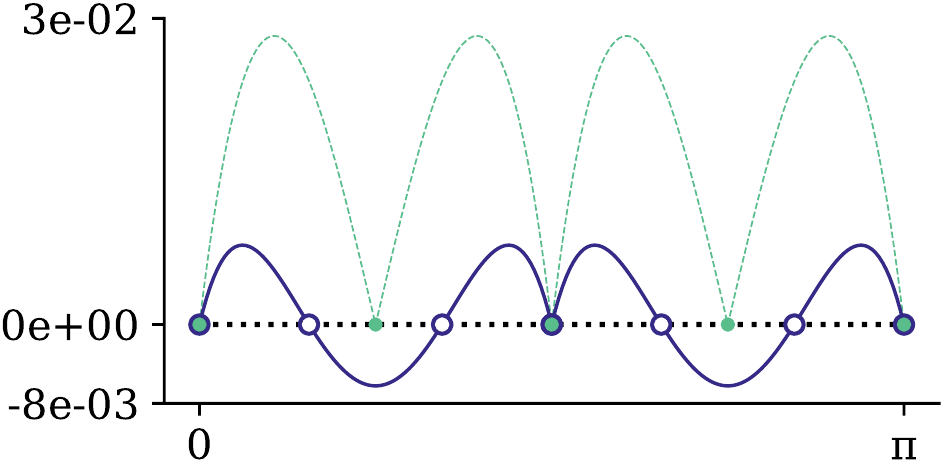} 
&
\includegraphics[
height=keepaspectratio,height=0.15\textwidth]{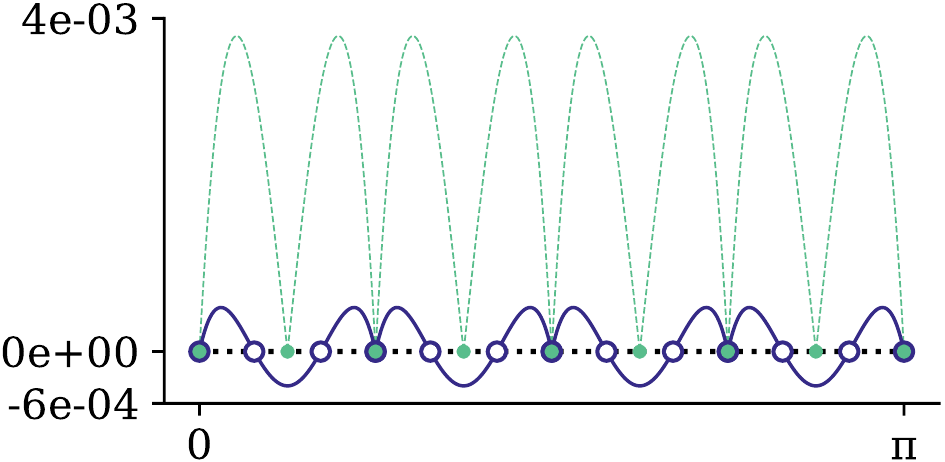} 
\\[10pt]

\end{tabular}

\caption{\label{fig:circle_1ele_fix_ref} 
approximating a semi-circle with the constrained disparity using $p=2$ meshes and R$=1,2, 4 $ elements, respectively. Top: optimized curves \hspace{2pt}   
\protect\tikz[baseline=-3pt]{\protect
\draw (0.0, 0.0)--(0.3, 0.0) [par_1, line width =1pt]} $\tilde \bx^\star$ \hspace{1pt}  and \hspace{1pt} 
\protect\tikz[baseline=-3pt]{\protect
\draw (0.0, 0.0)--(0.3, 0.0) [par_6, densely dashed, line width =1pt]}  $\ba\circ \tilde s^\star$. Bottom: error curves featuring internal roots
\protect\tikz[baseline=-3pt]{
\protect\fill(0.0, 0.0)[par_1] circle(2.5pt) ;
\protect\fill(0.0, 0.0)[white] circle(1.5pt) ;
}   \hspace{1pt}   
and interface points 
\protect\tikz[baseline=-3pt]{
\protect\fill(0.0, 0.0)[par_1] circle(2.0pt) ;
\protect\fill(0.0, 0.0)[par_6] circle(1.5pt) ;
}.}
\end{figure}

In Figure \ref{fig:circle_1ele_fix_ref} we show a semi-circle approximated with several elements and the corresponding point-wise errors before and after optimizing the constrained disparity. Note how although only internal nodes move during optimization, the error magnitude reduces. In Section \ref{sec:local_superconvergence}, we will see that the constrained disparity preserves the super-convergent rates from the original formulation. Moreover, fixing element interfaces transforms the optimization problem in R (total elements) independent copies, allowing optimization in parallel. At the end of this section, we discuss details on parallelization.

\subsection{Logarithmic barrier}
The commutative diagram shown in Figure \ref{fig:disp_diagram} holds if all mappings are diffeomorphisms. Since there are no constraints in formulation \eqref{eq:disp2}, in particular, diffeomorphism  $s$ is not actively enforced. The curve will tangle if elements in the parametric mesh $s$ are inverted. A possible solution is to add a constraint on $s'$ through the line search \cite{garimella2004, dovrev}. Alternatively, we can avoid curve tangling by introducing a log barrier function.

\textbf{Log barrier \cite[Ch.9] {nocedal}:} 
consider the nonlinear programming problem:
\[\min f(x) \quad 
\text{subject to} \quad c_i \geq 0,\quad i=1,\ldots,m.
\]
The logarithmic barrier is a  \emph{penalty} term that moves the optimizer away from violating the constraint $c_i$. For a given $\mu$, one can solve instead  
\[
\min\limits_x P(x;\mu) = \min\limits_x f(x) -  \mu \sum_{i=1}^m \log(c_i(x)). 
\]

\begin{figure}[h!]
\flushleft
\hspace{10pt}
\begin{tabular}{cc}
\begin{overpic}[
width = keepaspectratio, width =0.44\textwidth, abs,unit=1mm]
{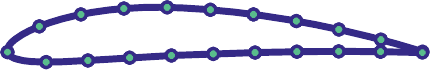} 
\put(50,10){\tikz{\draw (0.7,0) node {$\ba \circ s$}}}
\end{overpic}
&\hspace{10pt}
\begin{overpic}[
width = keepaspectratio, width = 0.44\textwidth]
{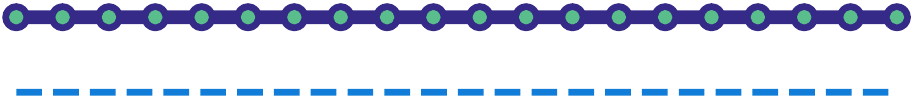} 
\put (95,12) {$s'$}
\end{overpic}
\end{tabular}
\caption{\label{fig:naca:log_barrier_init} initial parametrization of a NACA curve ($\ba\circ s$) and $s'$ values relative to the zero axis \raisebox{0.1cm}{\protect\tikz{\protect
\draw (0.0,0.0)--(0.4,0.0)[par_3,densely dashed, line width = 1.25pt]}}.}
\end{figure}

A valid curve parametrization follows the curve either forward or backwards. When we compose $\ba\circ  s$, we can ensure that the reparametrization preserves direction by fixing the sign of $s'$ (positive being forward and negative backwards) at the beginning of the optimization. For example, in Figure \ref{fig:naca:log_barrier_init} we show a NACA curve and the initial parametrization  $\ba\circ s = \ba$. Note that $s$ is not optimized yet, so it is simply a linear mapping (constant derivative). 

\vspace{20pt}

We are solving the following problem:  
\[
\min\limits_{\bx,s} E(\bx,s) \quad 
\text{subject to} \quad s'(\xi) > 0\quad \forall \xi.
\]
We need a continuous barrier function, so we introduce: 
\begin{equation}\label{eq:log_disp}
P(\bx,s;\mu) = E(\bx,s) - \mu\int\limits_{\mathcal M^R} \log (s'(\xi))d\xi.
\end{equation}
For each $\mu$, we optimize instead the following: 
\begin{equation}\label{eq:min_log_disp}
(\bx^\star,s^\star)  = \argmin\limits_{\bx,s} P(\bx,s;\mu),
\end{equation}
and use that if $s$ is not tangled, then
\[
\lim_{\mu\to 0} P(\bx,s;\mu) = E(\bx,s).
\]
\begin{remark}
In this case, $E$ stands for either the original (unconstrained)  $(\bx^\star,s^\star)$ or the constrained  ($\tilde \bx^\star,\tilde s^\star$) solutions. Both benefit from this technique. 
\end{remark}

\begin{figure}[h!]
\flushleft
\begin{tabular}{cc}
\textbf{Unconstrained} & 
\textbf{Log barrier}\\[8pt]
\begin{overpic}[
width = keepaspectratio, width =0.44\textwidth, abs,unit=1mm]
{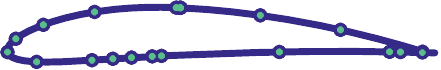} 
\end{overpic}
& \hspace{10pt}
\begin{overpic}[
width = keepaspectratio,  width =0.44\textwidth, abs,unit=1mm]
{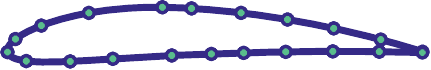} 
\put (45,8) {\begin{tikzpicture}
\draw (0,0)--(0.25,0) [par_1, line width = 1.5pt];
\fill (0,0) [par_1] circle (2pt);
\fill (0,0) [par_6] circle (1pt);

\draw (1.0,0.05) node {$\ba \circ s^\star$};
\end{tikzpicture}}
\end{overpic}
\\[15pt]
\raisebox{-18pt}{
\begin{overpic}[
height = 0.15\textwidth, width = 0.4\textwidth, abs,unit=1mm]
{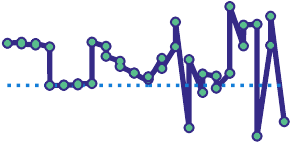} 
\end{overpic}}
& 
\begin{overpic}[
height = 0.05\textwidth, width = 0.4\textwidth, abs,unit=1mm]
{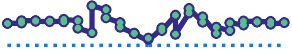} 
\put (50,5) {\begin{tikzpicture}
\draw (0,0)--(0.25,0) [par_1, line width = 1.5pt];
\fill (0,0) [par_1] circle (2pt);
\fill (0,0) [par_6] circle (1pt);
\draw (0.6,0.05) node {${s^\star}'$};
\end{tikzpicture}}
\end{overpic}
\end{tabular}
\caption{\label{fig:untangling} optimizing the NACA curve with (left) and without (right) the log barrier. The horizontal axis, \protect\tikz{\protect
\draw (0.0,0.0)--(0.5,0.0)[par_3,densely dashed, line width = 1.25pt]}, denotes $s'=0$. }
\end{figure}

In Figure \ref{fig:untangling} (left), we show the optimized NACA curve $\ba\circ s$ without a log barrier where the curve develops artificial loops. Observe how the derivative profile $s'$ crosses the zero axis in several locations. On the right, we show the same curve but optimized with the $\log$-barrier,  resulting in a valid curve. 

\vspace{5pt}

\begin{note}
In practice, the $\log$-barrier activates only if, at a particular Newton step, we detect a change in the sign of $s'$. This check is done by oversampling $s$ at each element. At that point, it retrieves the previous valid pair $(\bx_n,s_n)$ and solves instead the penalized problem $P(\bx,s;\mu_k),\ k = 1,\ldots,M$. 
\end{note}

\subsection{Globalized Newton with Zhang-Hager line search}
We find the solution to our minimization problem by combining Newton with a backtracking line search: 
\begin{align}
(\bx_{n+1},s_{n+1}) = (\bx_n,s_n) + \alpha_n \boldsymbol d_n.
\end{align}
Let $\nabla(\cdot)$ and $H(\cdot)$  denote the gradient and Hessian operators respectively and $\boldsymbol \delta_n$ a Newton step:
\begin{equation}\label{eq:newton_step}
H (E_n) \boldsymbol{\delta}_n = -\nabla(E_n).
\end{equation}
At each iteration, we ensure a descent direction using: 
\begin{align}\label{eq:descent}
\boldsymbol d_n :=\left\{
\begin{aligned}
\boldsymbol \delta_n,&\quad \text{if}\ \boldsymbol \delta_n \cdot \nabla E_n < 0,\\
- \Big[ diag(H(E_n)) \Big]^{-1} \nabla E_n,&\quad \text{otherwise}.
\end{aligned}\right.
\end{align}
As mentioned before, the usual line search choice is the \emph{Armijo} (monotone) rule:
 $\alpha \in(0,1]$ satisfies that
\begin{equation}
E_{n+1} < E_n + \alpha  10^{-4} \boldsymbol d_n \cdot \nabla E_n.
\end{equation}

Instead, we use a type of Zhang-Hager nonmonotone line search \cite{zhang}. Let  $C_0 = E_0,\ Q_0 = 1$ and define:
\begin{equation}\label{eq:cqn}
Q_{n+1} = \eta_n Q_n + 1\quad C_{n+1} = \frac{\eta_n Q_n C_n + E_{n+1} }{Q_{n+1}}.
\end{equation}
Notice that $\eta \equiv 0$ gives Armijo's monotone line search. On the other hand, $\eta \equiv 1$ gives an average-based rule:
\begin{equation}\label{eq:average_ls}
C_n = \frac{1}{n}\sum\limits_{i = 1}^n E_n, 
\end{equation}
which is the one that we will use in our experiments.
Finally,  $\alpha_n$ satisfies \textbf{Wolfe conditions:}
\begin{align}\label{eq:wolfe_nm}
E_{n + 1}  \leq C_n + \sigma_1 \alpha_n \nabla E_n \cdot d_n\\
\nabla E_{n+1} \cdot d_n \geq \sigma_2 \nabla E_n \cdot d_n
\end{align}
with $\sigma_1 = 10^{-4}$ and $\sigma_2 = 0.9$ as suggested in \cite{nocedal}.

\subsection{Implementation of the constrained disparity optimizer}

Now we provide implementation details for the proposed nonlinear optimization problem:
\begin{align}
\tilde E(\tilde \bx^\star,\tilde s^\star) = \min\limits_{\tilde \bx,\tilde s} || \tilde \bx -\ba \circ\tilde  s||_{\bs},
\end{align}
with $\tilde \bx$ and $\tilde s$ defined at each element (see equation \eqref{eq:const_meshes}) by: 
\begin{align}
\tilde \bx(\xi) = \bx^F(\xi)+\sum\limits_{i = 2}^{p} \tilde \bx_i N_i^p(\xi),\qquad
\tilde s(\xi) = s^F(\xi)+\sum\limits_{i = 2}^{q} \bx_i N_i^q(\xi).
\end{align} 
\begin{remark} The values $\bx^F$ and $s^F$ are fixed during optimization. 
\end{remark}

The main function in Algorithm~\ref{alg:cdm}, OPTIMIZE, takes several arguments. The first,  $\ba$, gives information about the curve. The initial meshes are given by $\bx_0$ and $s_0$, respectively. M is the number of outer iterations corresponding to the log barrier term $\mu$, and N is the maximum nonlinear iterations allowed. In our experiments, we use M=6 with $\mu$ decreasing by a factor of $10^{-2}$ at each step. The function returns the optimal pair ($\tilde \bx_n,\ \tilde s_n)$.

The optimization routine proceeds as follows. First, we initialize the log barrier variables (line 2). The outer loop (starting at line 3) corresponds to the penalized problem (equation \eqref{eq:log_disp}). The inner loop (lines 5-23) is the backtracking Newton scheme minimizing $\tilde E(\tilde \bx^\star,\tilde s^\star)$. At each iteration, we compute the gradient and Hessian (line 9) and check the stopping criteria (lines 7-8). In our experiments, it corresponds to $tol = 10^{-12}$. Then, using a descent direction (lines 9-11), we loop until the line search condition is met (lines 13-16). Finally, we sample $s'$ for any changes in its sign and activate the log barrier if necessary (lines 17-21). 

\algrenewcommand\algorithmicindent{2em}%
\begin{algorithm}[h!]

\caption{\label{alg:cdm}Constrained disparity minimization }
\begin{algorithmic}[1]

\Function{Optimize}{$\ba$, $\bx_0$, $s_0$, $p$, $q$, M, N, tol}
    \State $\mu = 0$; logBarrier = $false$;
     \For{$m=1:$M}
     	 \State$ \mu = \mu \cdot 10^{-2}$
        \For{n=1:N}
            \State $\nabla E_n, H (E_n) \gets $\Call{GradHess}{$\ba,\ \bx_{n},\ s_{n}$, $\mu_k$}
            \If{$|\nabla (E_n)| \ <\ tol$}\State \textbf{break}
            \EndIf
            \State $\boldsymbol \delta = -H(E_n)^{-1}\nabla E_n$ 
            \If{$\boldsymbol \delta \cdot\nabla E_n \ > 0$}
            \State$ \boldsymbol \delta = - \Big[ diag(H(E_n)) \Big]^{-1} \nabla E_n$
            \EndIf
            \State $\beta = 1$
            \Repeat
                \State $(\tilde x_{n+1},\ \tilde s_{n+1}) = (\tilde x_n,\ \tilde s_n) + \beta \boldsymbol{\delta} $
               \State  $\beta = \beta/2$
            \Until{\Call{ZhangHager}{$\tilde x_{n+1},\ \tilde s_{n+1}$} = true}
             \If{\Call{any}{$(sign(\tilde s'_{n+1})\neq sign(\tilde s'_0))$}}
                \State $(x_{n+1},\ s_{n+1}) = (x_n,\ s_n)$
                \State$\mu = || \tilde \bx_n - \ba \circ \tilde s_n||^2_{\bs}$
                \State logBarrier = $true$
                \State \textbf{break}
            \EndIf
         \EndFor
         \If{logBarrier = $false$}\State \textbf{break}
        \EndIf
    \EndFor
    
    \State \Return {$\tilde\bx_{n},\tilde s_n$}
\EndFunction
\end{algorithmic}
\end{algorithm}

\subsection{Parellelization}
Before we discuss the optimizatino in parallel, we describe our software environment. The optimization requires extracting curves from the models and querying information about the derivatives. We achieve this through the EGADs open-source geometry kernel \cite{egads} and EGADSlite \cite{egadslite}, designed to handle geometry in a parallel environment. We have developed a \emph{Julia} interface that is now available through the ESP distribution \cite{esp}. Further, we use \emph{Julia}'s distributed memory packages combined with the function $\verb|@spawnat|$ at the for loops. In our experience, this approach outperforms the $\verb|Threads.@thread|$ multi-threading utility.

CAD models generally consist of hundreds of curves. A natural parallelization would be distributing the curves among the available CPUs as described in Algorithm \ref{alg:par_curves}. First, we load the CAD data and extract the curves from the model (lines 2-3). Then, we distribute the curves among the available CPUs (lines 4-7). Similarly, we distribute arrays corresponding to the physical and parametric meshes (lines 8-9). Finally,  we optimize each curve in a parallel loop (lines 11-12).

\begin{algorithm}
\caption{\label{alg:par_curves}Parallel optimization by curves}

\begin{algorithmic}[1]
\Function{CurveParallel}{CADfile, $p$, $q$, $n$, M, N}
	\State $\text{aux} \gets$ \Call{load file}{CADfile}
	\State $curves  \gets$ \Call{findall}{aux != straight line} \Comment{exclude linear meshes}
	\State $\ell \gets $\Call{length}{curves} 
	\State $tasks \gets\min(nworkers(), \ell)$ \Comment{nworkers() = CPUs available}
	\State $workers = [1:tasks;]$
	\State $ids     = distribute([1:\ell;],\ procs = workers)$
	\State $\vec \bx    = distribute(\verb|Array{3Dmesh(n)}(\ell)|,\ procs = workers)$
	\State $\vec s      = distribute(\verb|Array{1Dmesh(n)}(\ell)|,\ procs = workers)$
	\State $\verb|const|\ lc = localpart$ \Comment{alias for localpart}
    \For{j = 1 : nw}
		\State $\verb|@spawnat|\ j\quad$\Call{Optimize}{$curves(lc(ids))$, $lc(\vec \bx)$, $lc(\vec s)$, M, N}
    \EndFor
    \State \Return {$\vec \bx, \vec s$}
\EndFunction
\end{algorithmic}
\end{algorithm}
\begin{note}
Algorithm \ref{alg:par_curves} can be used on the unconstrained disparity as well, changing the function "optimize" so that the element interfaces can move. This is the only possible parallelization for such formulation since the optimization requires solving all elements at once.
\end{note}

We have discussed how the constrained disparity transforms the optimization problem into independent $R$ copies of a lower dimensional problem. Here $R$ is the number of elements, and dimension refers to the size of the  optimization problem. Assuming a suitable initial element partition per curve, we can distribute and optimize the elements among the available cores. The routine is similar to Algorithm \ref{alg:par_curves} and is described in Algorithm \ref{alg:par_elements}. In this case, we distribute the mesh elements among the available workers (lines 3-6) and then perform the optimization loop in parallel (lines 7-10).

\begin{algorithm}

\caption{\label{alg:par_elements}Parallel optimization by elements}

\begin{algorithmic}[1]

\Function{ElementParallel}{$\ba$, $\bx$, $s$, M, N}
	\State $n \gets $\Call{NumberOfElements}{$\bx$}
	\State $tasks \gets\min(nworkers(), n)$ \Comment{nworkers() = CPUs available}
	\State $workers = [1:tasks;]$
	\State $\bx = distribute(\verb|Array{3Dnodes}(n)|,\ procs = workers)$
	\State $s   = distribute(\verb|Array{1Dnodes}(n)|,\ procs = workers)$
	\State $\verb|const|\ lc = localpart$ \Comment{alias for localpart}
    \For{j = 1 : nw}
		\State $\verb|@spawnat|\ j\quad$\Call{Optimize}{$\ba$, $lc(\bx)$, $lc(s)$, M, N}
    \EndFor
    \State \Return {$\bx,s$}
\EndFunction
\end{algorithmic}
\end{algorithm}

\section{Constrained disparity: exploiting local higher accuracy}
 \label{sec:local_superconvergence}
In this section, we discuss the performance of the constrained disparity regarding accuracy. We discuss the super-convergence phenomena starting with a global study of the constrained disparity. Then, we study the local error focusing on a single element and show that we can attain super-convergence by optimizing only the internal nodes.

\begin{figure}[h!]
\flushleft
\begin{tabular}{cc}

\hspace{5pt}
\textbf{Circle - 2D}  
\raisebox{-10pt}{\includegraphics
[height=keepaspectratio,height=0.08\textwidth]
{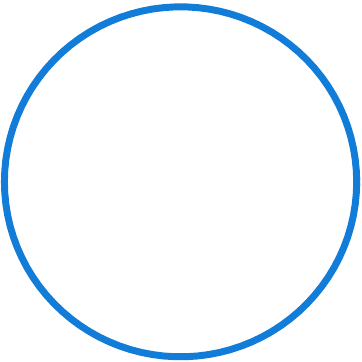}}
&
\hspace{20pt}
\textbf{Sphere arc- 3D} 
\raisebox{-10pt}{\includegraphics
[height=keepaspectratio,height=0.09\textwidth]
{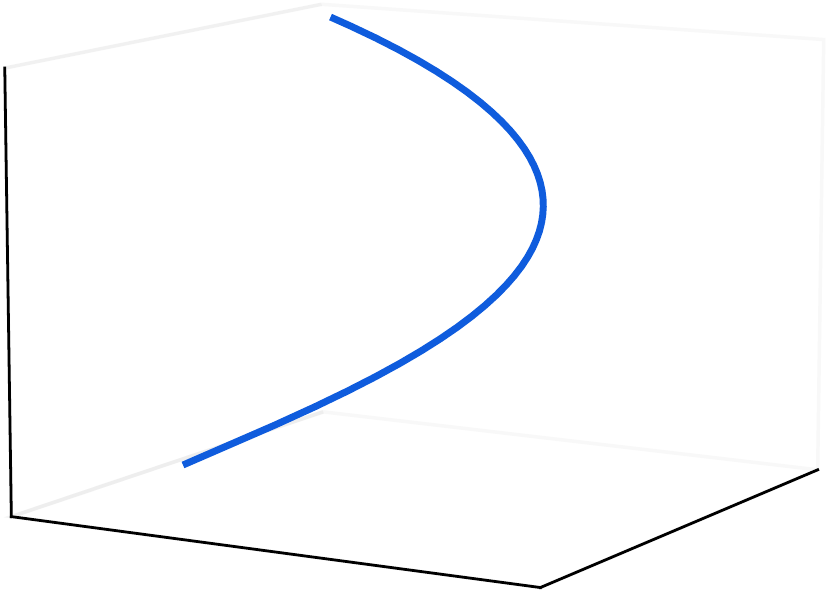} }
\\[20pt]
\hspace{5pt}
\begin{overpic}[
width=keepaspectratio,width=0.31\textwidth]
{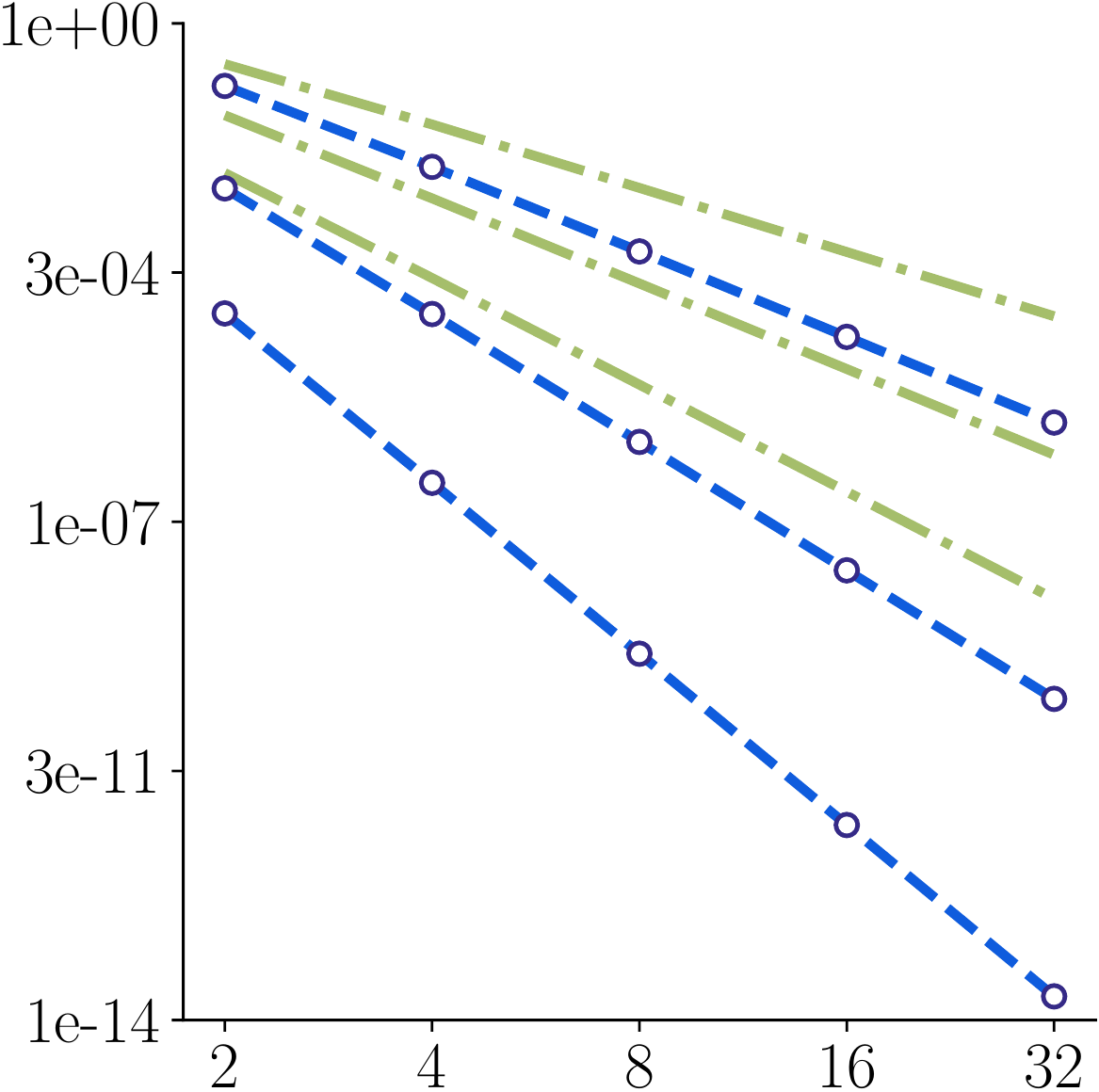}

\put(66,35){\tikz{\draw (0,0.2)--(0.4,-0.1)--(0.4,0.2)--cycle [par_2, line width = 1.pt]}}
\put(63,31){\tikz{\draw (0,0.2) node{\textcolor{par_2}{$\scriptstyle{\boldsymbol{8}}$}}}}
\put(66,55){
\tikz{\draw (0,0.2)--(0.4,-0.05)--(0.4,0.2)--cycle [par_2, line width = 1.pt]}}
\put(63,50){
\tikz{\draw (0,0.2) node{\textcolor{par_2}{$\scriptstyle{\boldsymbol{6}}$}}}
}

\put(66,77){
\tikz{\draw (0,0.2)--(0.4,0.05)--(0.4,0.2)--cycle [par_2, line width = 1.pt]}}
\put(63,70){
\tikz{\draw (0,0.2) node{\textcolor{par_2}{$\scriptstyle{\boldsymbol{4}}$}}}
}
\put(100,59){\scriptsize{$p$=2}}
\put(100,34){\scriptsize{$p$=3}}
\put(100,9){\scriptsize{$p$=4}}
\put(-10,45){\rotatebox{90}{\scriptsize{$||\cdot||_{\bs}$}}}
\put(45,-8){\scriptsize{Elements}}
\end{overpic}
& \hspace{20pt}
\begin{overpic}[
width=keepaspectratio,width=0.3\textwidth]
{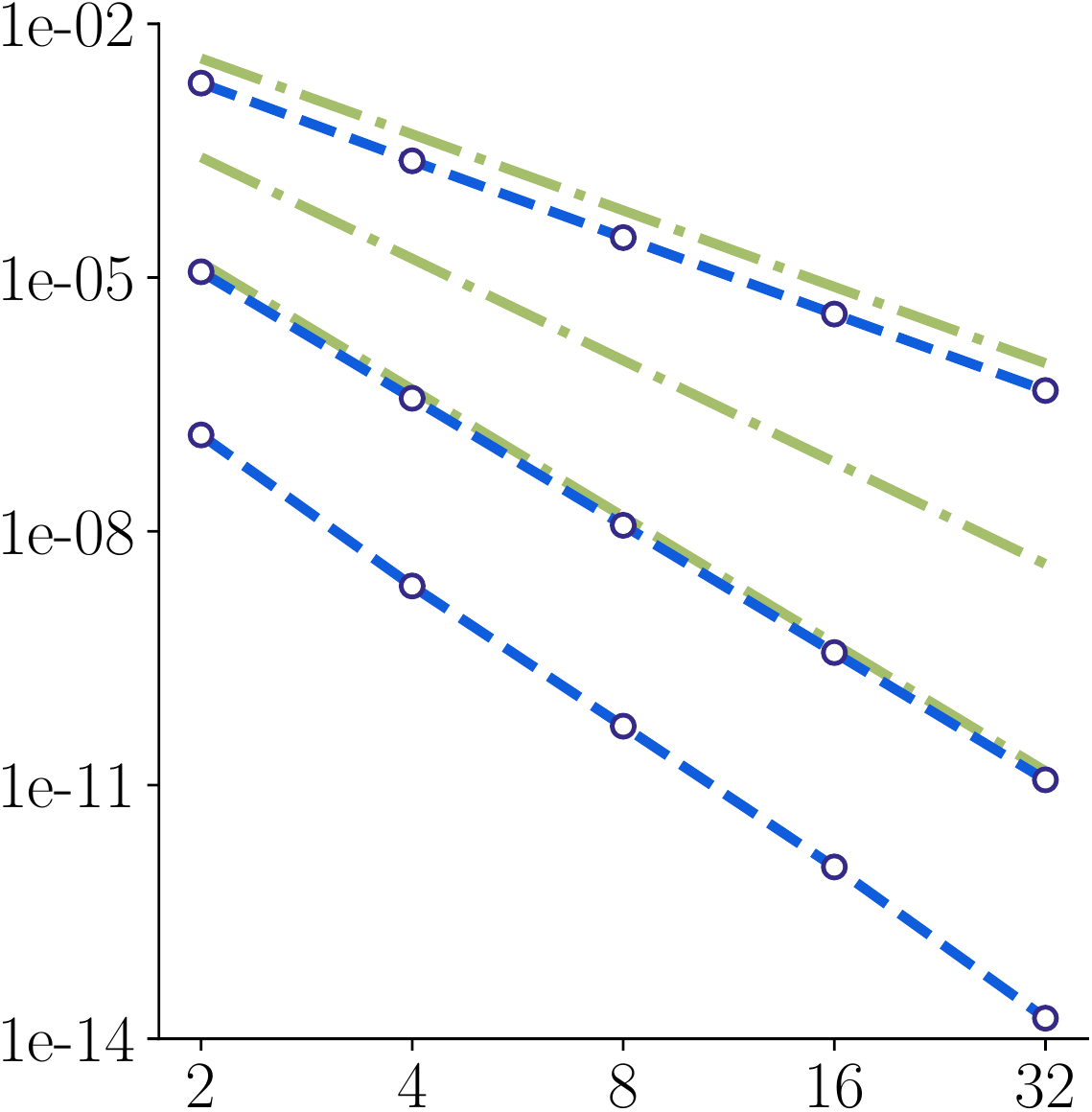}
\put(95,82){
\tikz{\draw (0,0.)--(0.5,0.0) [par_7, densely dashdotted, line width = 1.5pt]}}
\put(105,75){
\tikz{\draw (0.8,0.0) node{$|| \bx - \ba||_\sigma$}}}
\put(95,100){
\tikz{\draw (0,0.)--(0.5,0.0) [par_2, densely dashed, line width = 1.5pt]}}
\put(105,92){
\tikz{\draw (0.8,0.0) node{$||\tilde \bx^\star - \ba \circ \tilde s^\star||_\sigma$}}}

\put(66,28){
\tikz{\draw (0,0.2)--(0.4,-0.08)--(0.4,0.2)--cycle [par_2, line width = 1.pt]}}
\put(63,24){
\tikz{\draw (0,0.2) node{\textcolor{par_2}{$\scriptstyle{\boldsymbol{6}}$}}}
}

\put(66,48){
\tikz{\draw (0,0.2)--(0.4,-0.04)--(0.4,0.2)--cycle [par_2, line width = 1.pt]}}
\put(63,43){
\tikz{\draw (0,0.2) node{\textcolor{par_2}{$\scriptstyle{\boldsymbol{5}}$}}}
}

\put(66,78){
\tikz{\draw (0,0.2)--(0.4,0.05)--(0.4,0.2)--cycle [par_2, line width = 1.pt]}}
\put(63,72){
\tikz{\draw (0,0.2) node{\textcolor{par_2}{$\scriptstyle{\boldsymbol{4}}$}}}
}
\put(100,63){\scriptsize{$p$=2}}
\put(100,28){\scriptsize{$p$=3}}
\put(100,8){\scriptsize{$p$=4}}
\put(45,-8){\scriptsize{Elements}}
\end{overpic}
\end{tabular}\\[5pt]
\caption{\label{fig:constrained_rates_fix} slopes (log-log) of the constrained optimized disparity $||\cdot ||_{\bs}$, for several mesh refinements showing the super-convergent rates  \protect\tikz{\protect\draw(0.0,0.3)--(0.4,0.1)--(0.4,0.3)--cycle[par_2, line width = 1.pt]} :
$2p$ in 2D, and $\lfloor\frac 32 (p-1)\rfloor$ in 3D.}
\end{figure}

\begin{figure}[h!]
\flushleft
\vspace{10pt}
\hspace{90pt}
\begin{overpic}
[width = keepaspectratio, width = 0.35\textwidth]{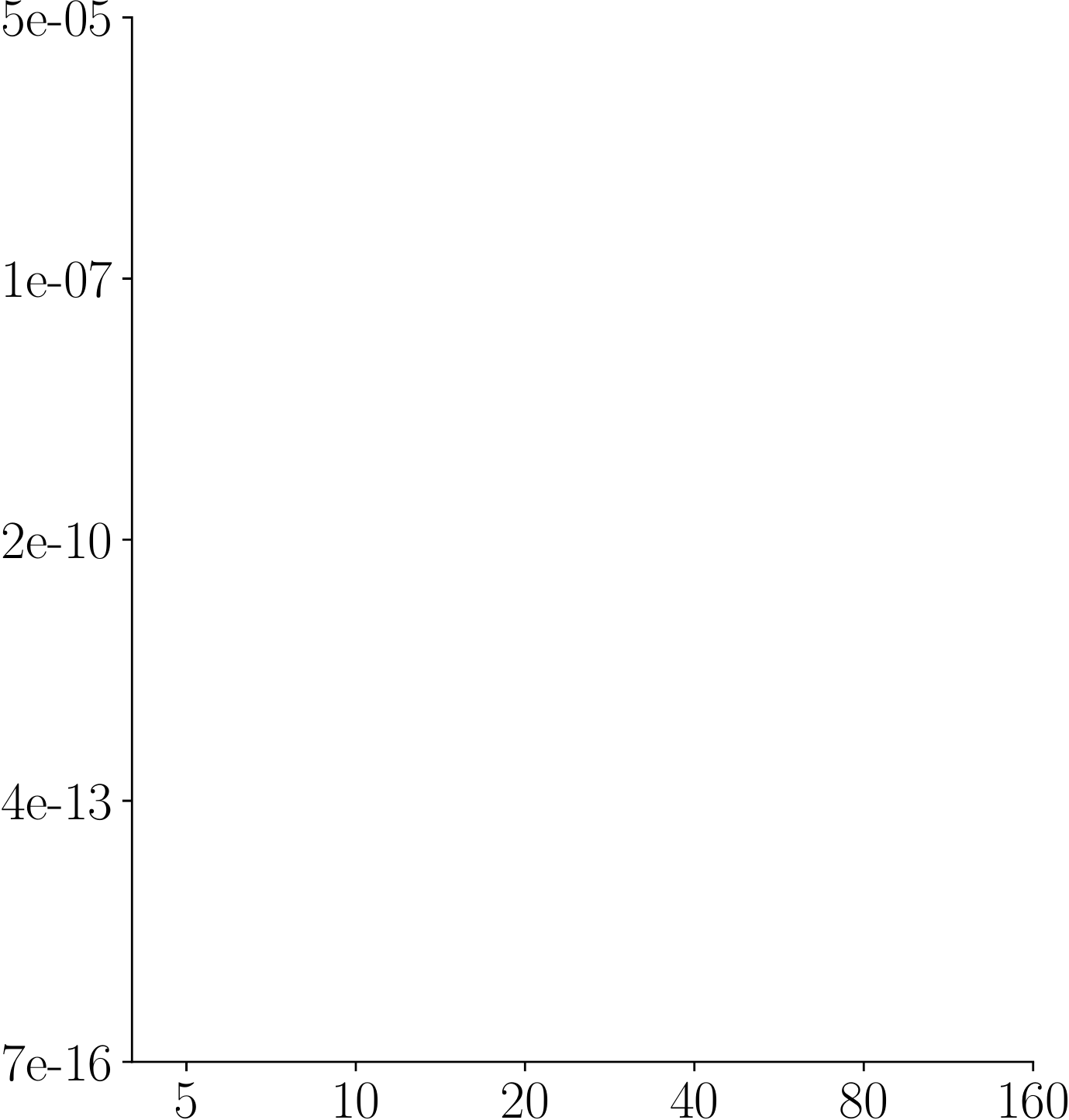}
\put(-60,80){\includegraphics[width = keepaspectratio, width = 0.2\textwidth]
{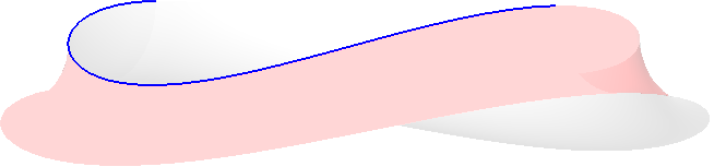}}
\put(9.6,5.1){
\includegraphics[width = keepaspectratio, width = 0.32\textwidth, abs]{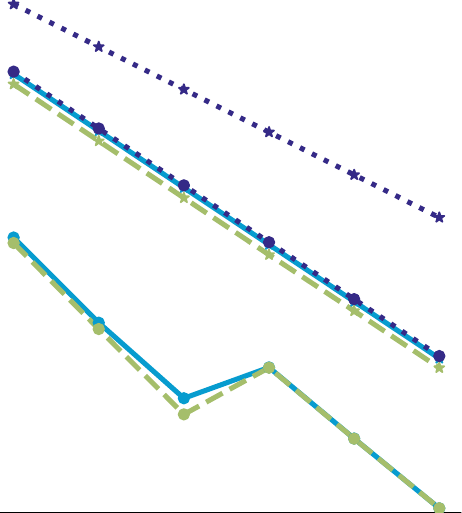}}
\put(-15,40){
\rotatebox{90}{$\boldsymbol{||\cdot ||_{\bs}}$}}

\put(97,58){\textcolor{par_1}{$\scriptstyle{\boldsymbol{p=2}}$}}

\put(97,36){\textcolor{par_1}{$\scriptstyle{\boldsymbol{p=3}}$}}

\put(97,30){\textcolor{par_3}{$\scriptstyle{\boldsymbol{p=2}}$}}

\put(97,5){\textcolor{par_3}{$\scriptstyle{\boldsymbol{p=3}}$}}
\put(42,-8){Elements}
\put(120,80){
\begin{tikzpicture}
\draw (0,0)--(0.4,0) [par_1, densely dotted, line width = 1pt];
\draw (1.3,0) node{$||\bx-\ba||_{\bs}$};
\end{tikzpicture}
}
\put(120,65){
\begin{tikzpicture}
\draw (0,0)--(0.4,0) [par_4, line width = 1pt];];
\draw (1.8,0) node{$||\tilde\bx^\star-\ba\circ \tilde s^\star||_{\bs}$};
\end{tikzpicture}
}
\put(120,50){
\begin{tikzpicture}
\draw (0,0)--(0.4,0) [par_7, densely dashed, line width = 1pt];];
\draw (1.8,0) node{$||\bx^\star-\ba\circ  s^\star||_{\bs}$};
\end{tikzpicture}
}
\end{overpic}
\vspace{10pt}
\caption{\label{fig:lazo_rates} 
convergence plots showing the slopes (log-log) of the $||\cdot ||_{\bs}$ norm for the top curve of the CAD model (marked in blue) using direct interpolation (dotted dark blue) vs. optimizing the constrained (solid light blue) and the original (dashed green) disparities.}
\end{figure}
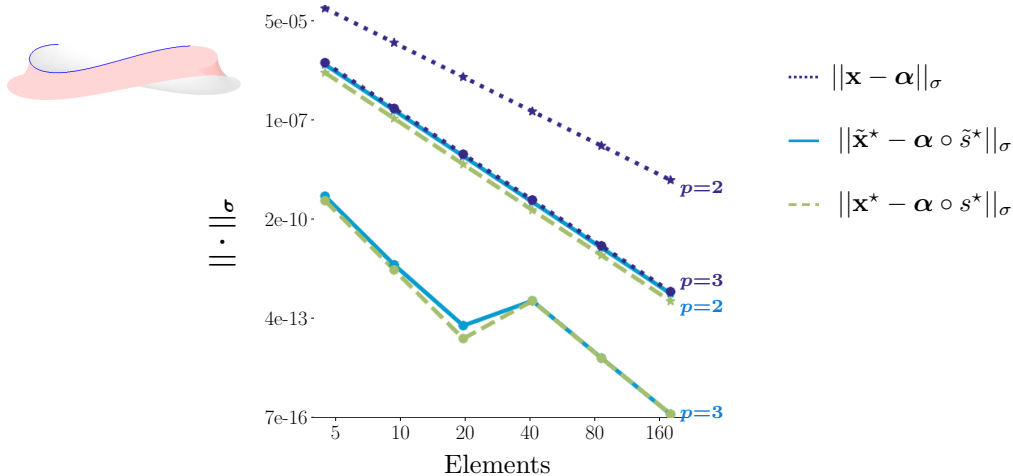

In Figure \ref{fig:constrained_rates_fix}, we show convergence plots for a circle and a sphere arc for $p=2,3,4$ and five mesh refinements. Notice that for the 2D case, as for the original disparity (Figure \ref{fig:circle_rates_optimal}), we raise the order from the expected $p+1$ to $2p$. For the 3D case, we obtain $\lfloor\frac 32(p-1)\rfloor + 2$ order. In Figure \ref{fig:lazo_rates} we perform a similar study using a curve from a CAD model as the target geometry. Although the constrained disparity is slightly larger, the order of accuracy is the same as for the original (unconstrained) problem, and both significantly improve the initial approximation.

\subsection{Planar curves: local error for a single element}
\begin{figure}[h]

\flushleft
\hspace{55pt}
\begin{tabular}{ccc}
\multicolumn{3}{c}{
\begin{tikzpicture}
\draw (0.0, 0.0)--(0.45, 0.0) [par_3, densely dashed, line width = 0.75pt];
\draw (0.7,0.0) node {$\be$ };

\draw (2.0, 0.0)--(2.45, 0.0) [par_5, densely dashed, line width = 0.75pt];
\draw (2.8,0.05) node {$\be^\star$ };

\draw (4.0, 0.0)--(4.45, 0.0) [par_1, line width =0.75pt];
\draw (4.8,0.05) node {$\tilde \be^\star$ };
\end{tikzpicture}}
\\
\begin{overpic}[
width=keepaspectratio,width=0.24\textwidth]{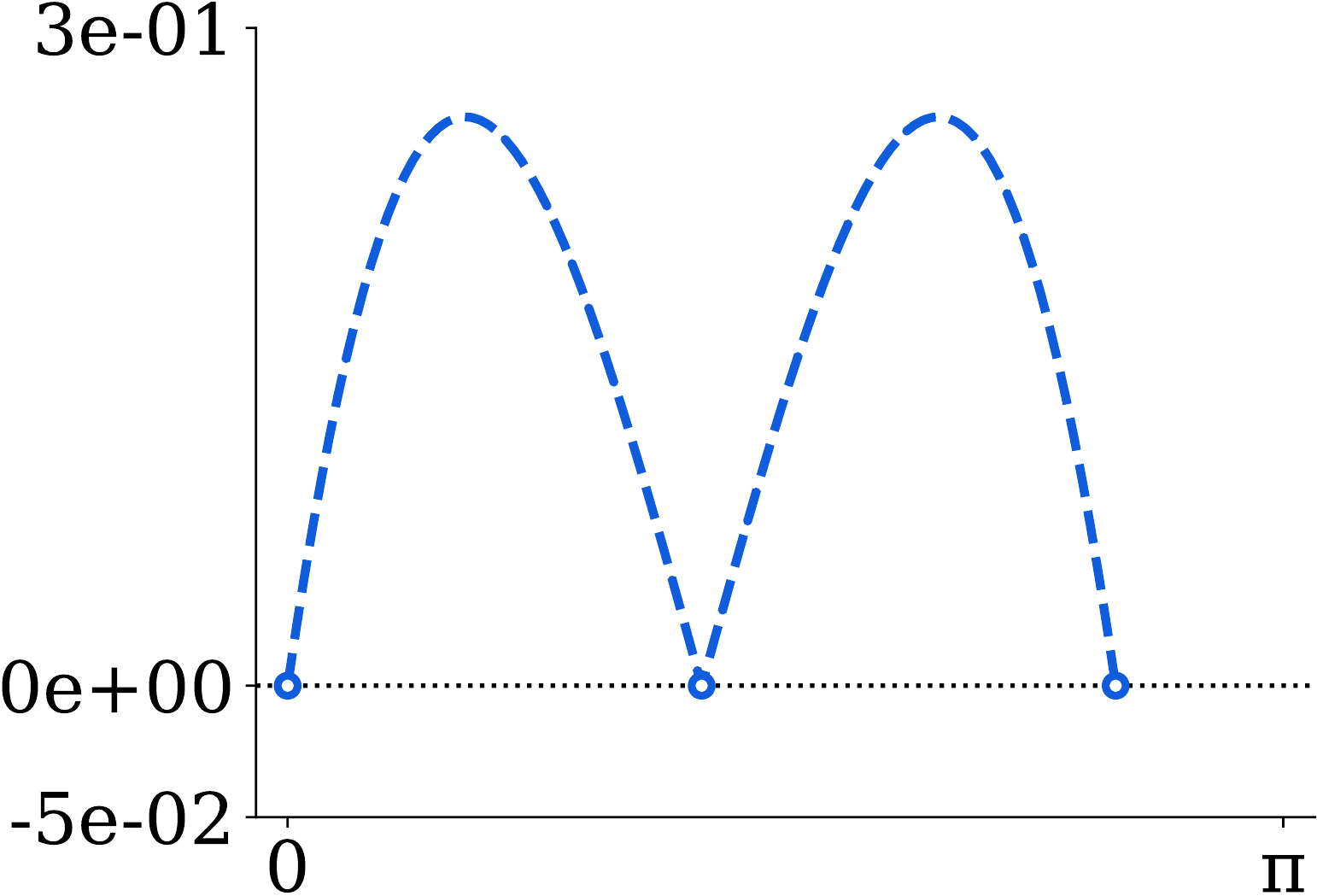} 
\put(-70,40){
\textbf{Initial}}
\end{overpic}
&
\includegraphics[
width=keepaspectratio,width=0.24\textwidth]{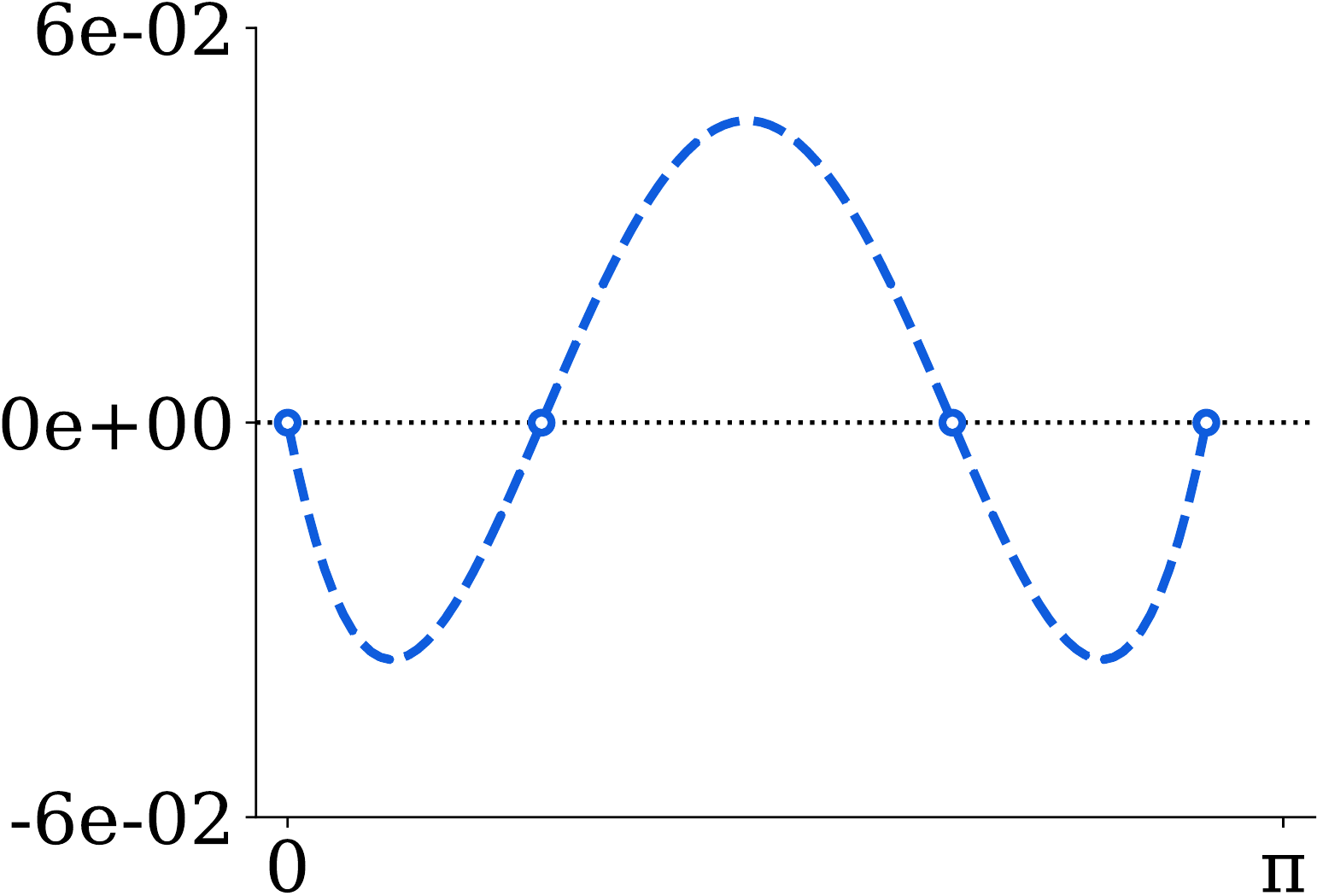} 
& 
\includegraphics[
width=keepaspectratio,width=0.24\textwidth]{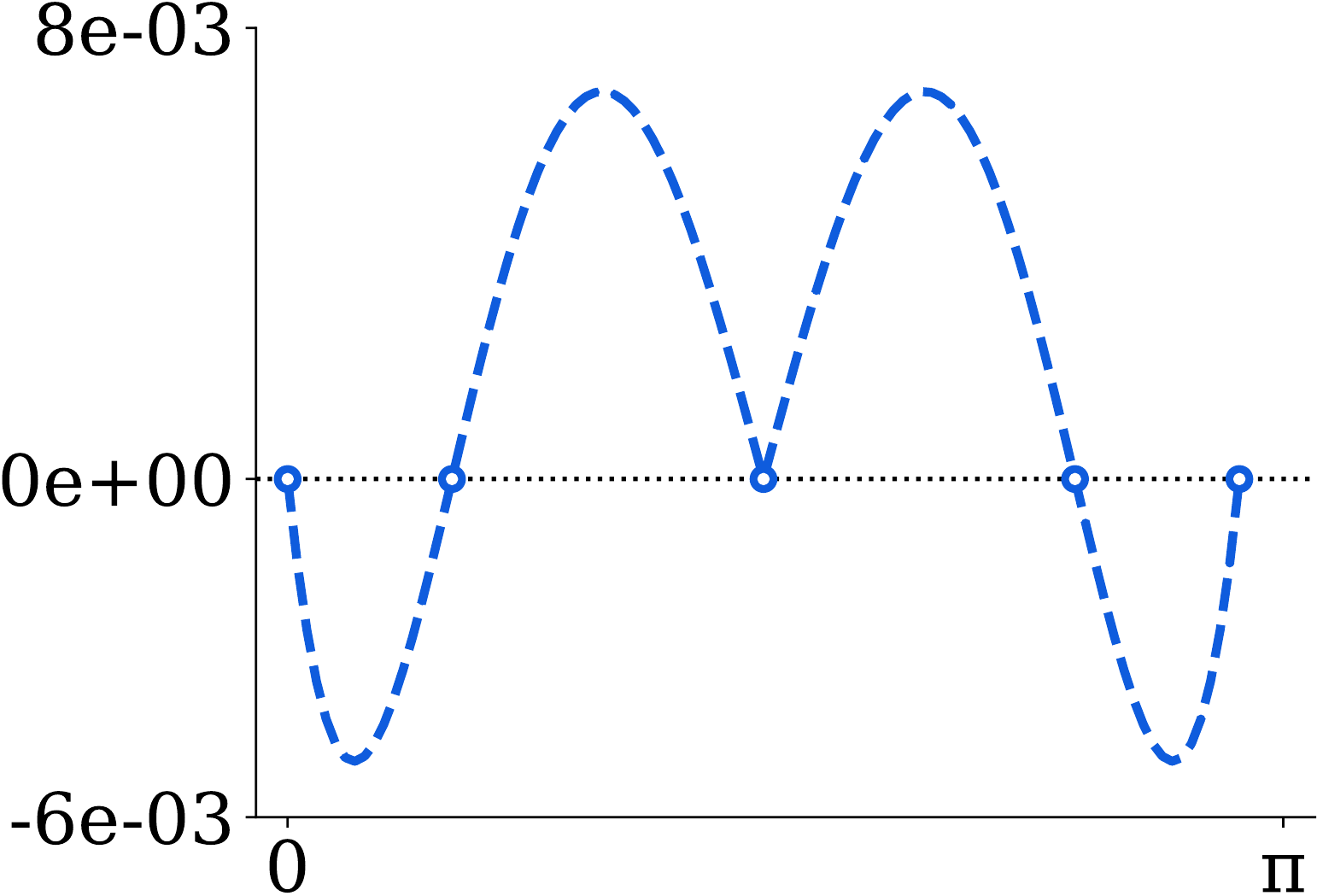}
\\
\begin{overpic}[
width=keepaspectratio,width=0.25\textwidth]{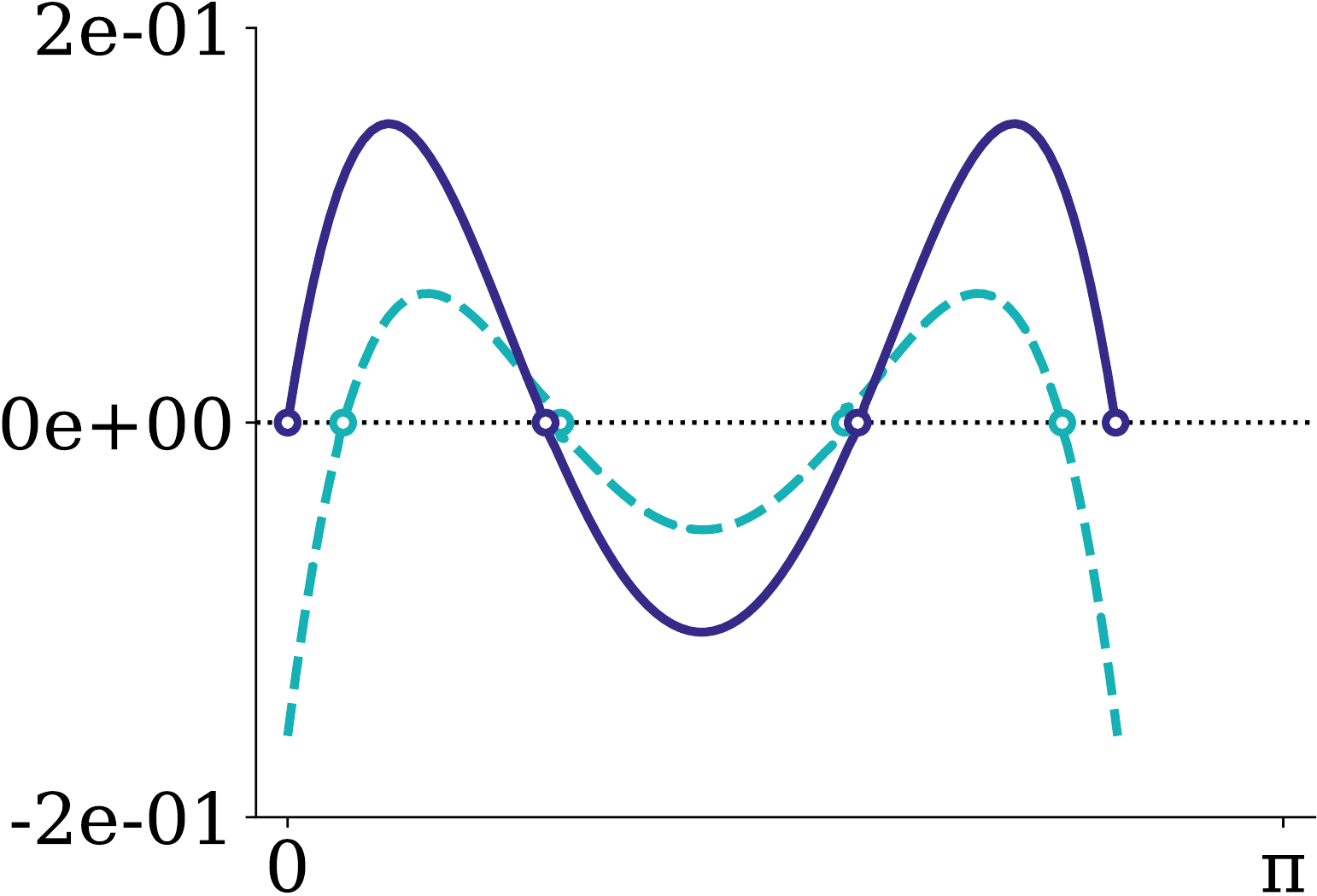} 
\put(-72,40){
\textbf{\small{Optimized}}}

\put(40,-10){$  \boldsymbol{ \scriptsize{p=2}}$}
\put(155,-10){$ \boldsymbol{ \scriptsize{p=3}}$}
\put(265,-10){$ \boldsymbol{ \scriptsize{p=4}}$}
\end{overpic}
&
\includegraphics[
width=keepaspectratio,width=0.25\textwidth]{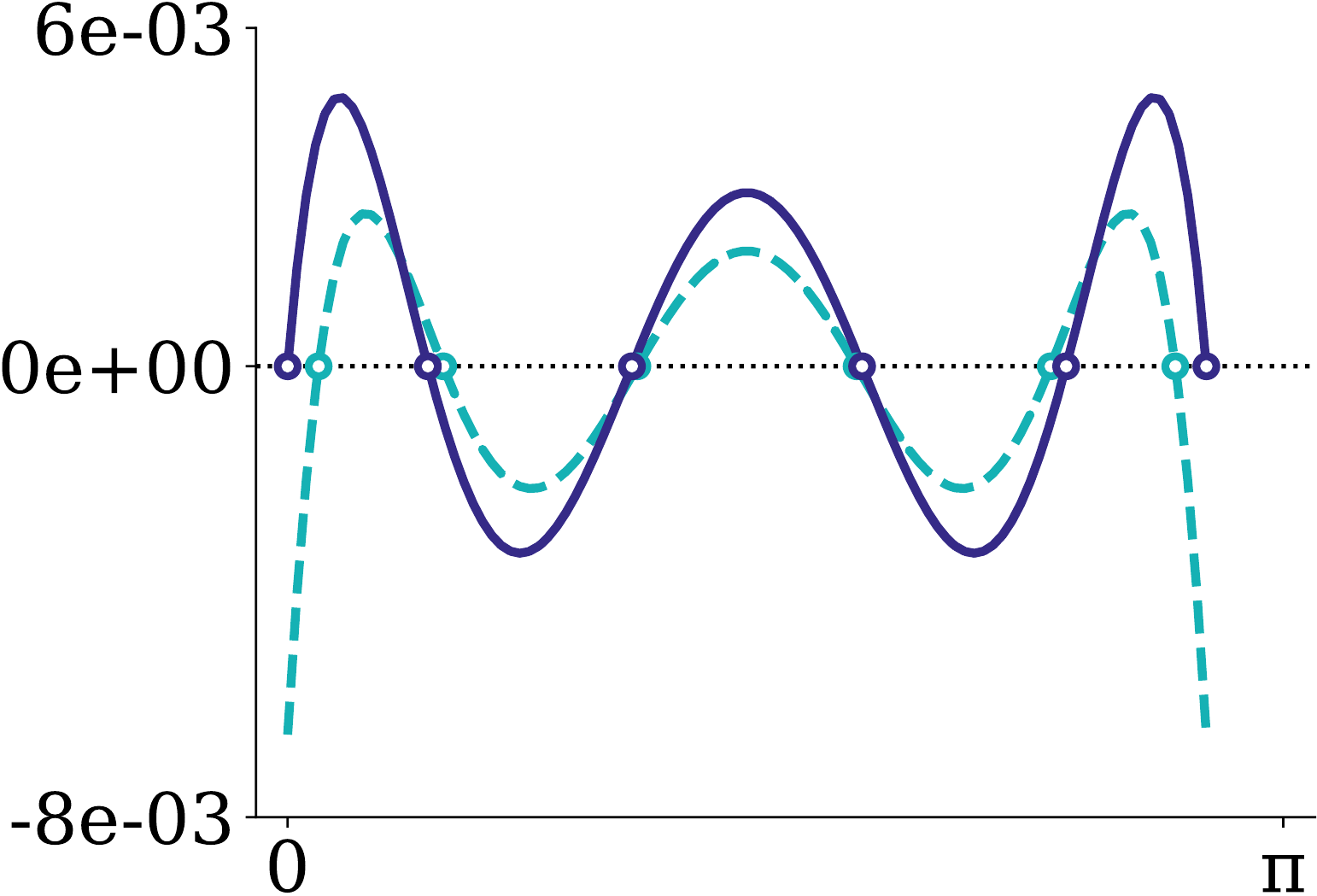} 

&
\includegraphics[
width=keepaspectratio,width=0.25
\textwidth]{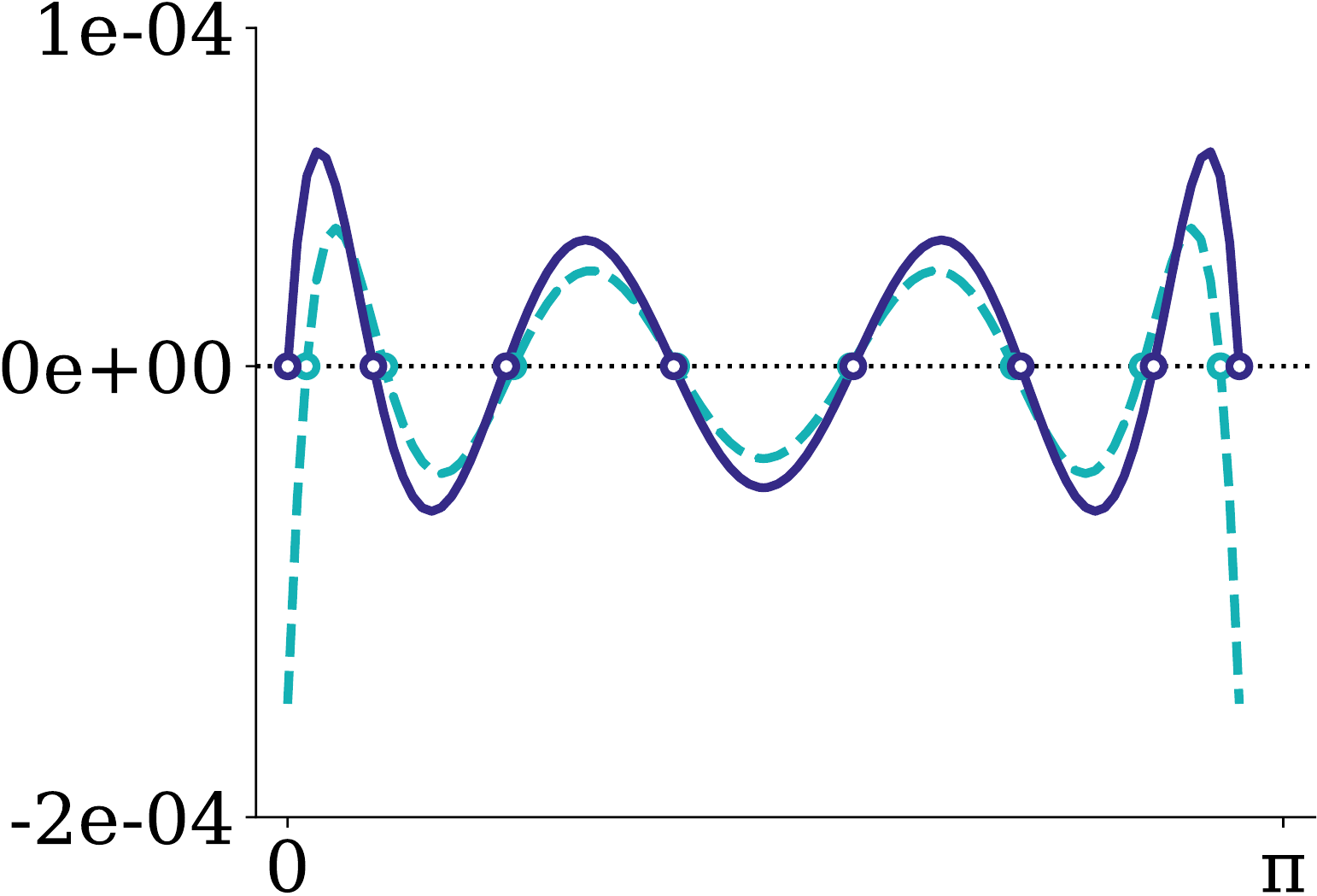}\\[4pt]
\end{tabular}

\caption{\label{fig:circle_1ele_roots} point-wise error plots, $\be = |\bx-\ba \circ s|$ (resp. $\cdot ^\star,\ \tilde \cdot ^\star$), approximating a semi-circle with a single element before (top) and after optimizing (bottom) the constrained (
\protect\tikz[baseline=-3pt]{\protect
\draw (4.0, 0.0)--(4.45, 0.0) [par_1, line width =1pt]}
$\ \tilde\be^\star$) and the original disparity (
\protect\tikz[baseline=-3pt]{\protect\draw
[par_5, densely dashed, line width = 1pt]  (2.0, 0.0)--(2.45, 0.0)}
 $\ \be^\star$). The $y$-axis denotes the magnitude of the error.}
\end{figure}

Here, we focus on a single element and study the local behavior of the optimizer. We use a semi-circle as the target curve to make the plots clearer. In Figure \ref{fig:circle_1ele_roots} we show point-wise error plots for several polynomial degrees. The initial approximation (top) is an interpolating polynomial of degree $p$, and the error curve has the expected behaviour: $p+1$ roots. The bottom plots show the results optimizing with fix and free interfaces. Notice that although the fixing interfaces produces slight larger errors, both solutions behave similarly: the curves have $2p$ roots (instead of $p+1$).

Now let us discuss the curve in terms of the Frenet frame. Denote  $\{\boldsymbol t,\boldsymbol n\}$ the curve tangent and normal vectors, respectively.
For planar curves, we can decompose the error $\be = \bx - \ba \circ s $ along these directions:  
\[\be = (\be\cdot \boldsymbol t)\boldsymbol t + (\be\cdot \boldsymbol n)\boldsymbol n.\]

For a single element, since the parametric mesh $s$ uses polynomials of degree $q$, we have a total of $q+1$ degrees of freedom. On the other hand, the physical mesh uses polynomials of degree $p$ in $\mathbb R^2$, so it has
$2(p + 1)$ degrees of freedom. Fixing the interfaces gives a total of $q+1-2 + 2(p + 1-2)$. Hence, we have $(q-1) + (2p-2)$ equations optimizing the disparity. 

\begin{figure}[h!]
\flushleft
\begin{tabular}{ccc} 

\mc{3}{
\hspace{40pt}\begin{tikzpicture}
\draw (0.0,0.0)--(0.5,0.0) [par_1, line width=1pt];
\draw (2.0,0.0)--(2.5,0.0) [densely dashed, par_4,line width=1pt];
\draw (1.0,0.0)  node {$\boldsymbol t$};
\draw (3.0,0.0)  node {$\boldsymbol n$};
\end{tikzpicture}} \\[5pt]

\hspace{40pt}
\begin{overpic}[
width=keepaspectratio,width=0.26\textwidth]{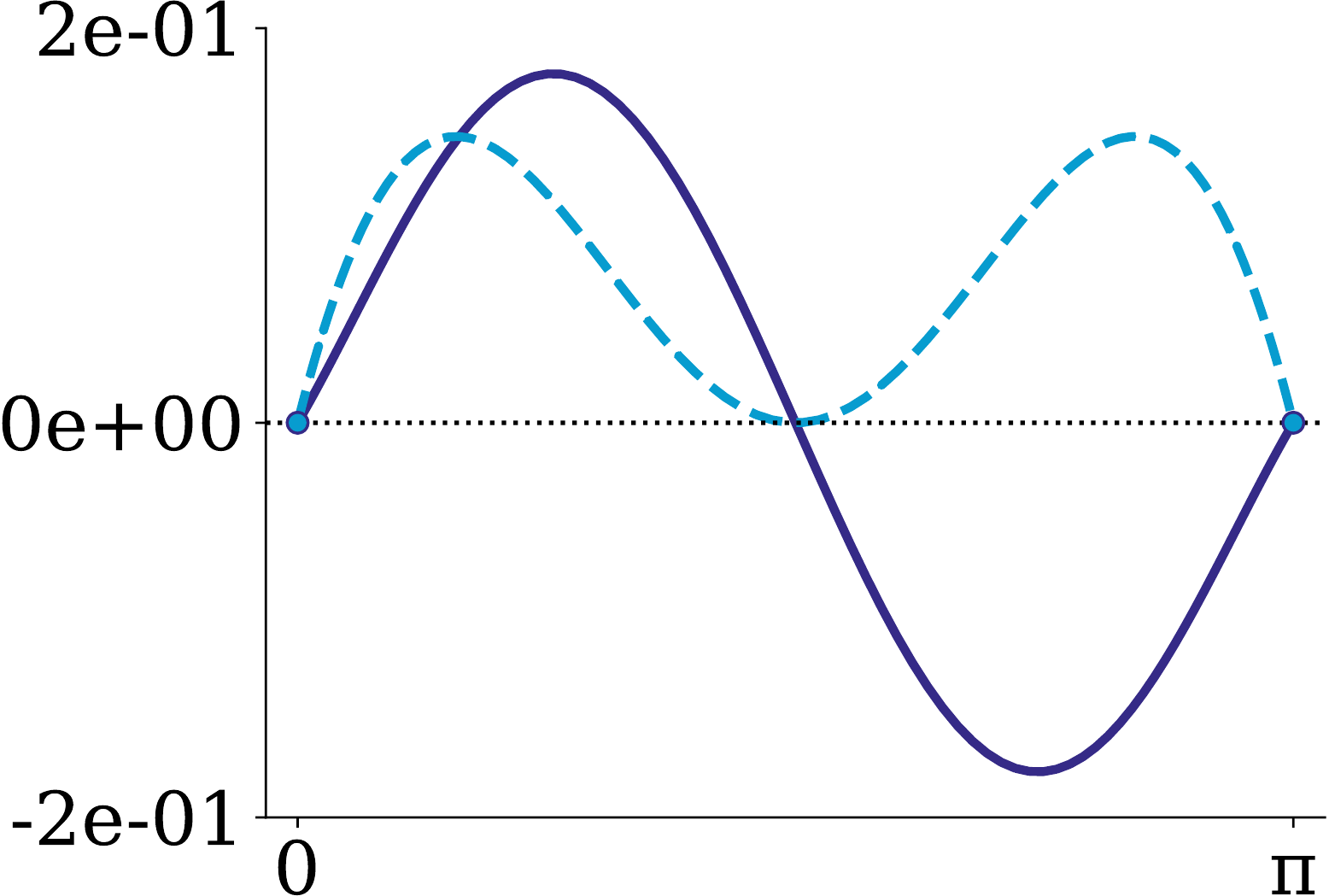} 
\put(-45,30){\fbox{$\boldsymbol{p = 2}$}}
\end{overpic} 
& 
\includegraphics[width=keepaspectratio,width=0.26\textwidth]{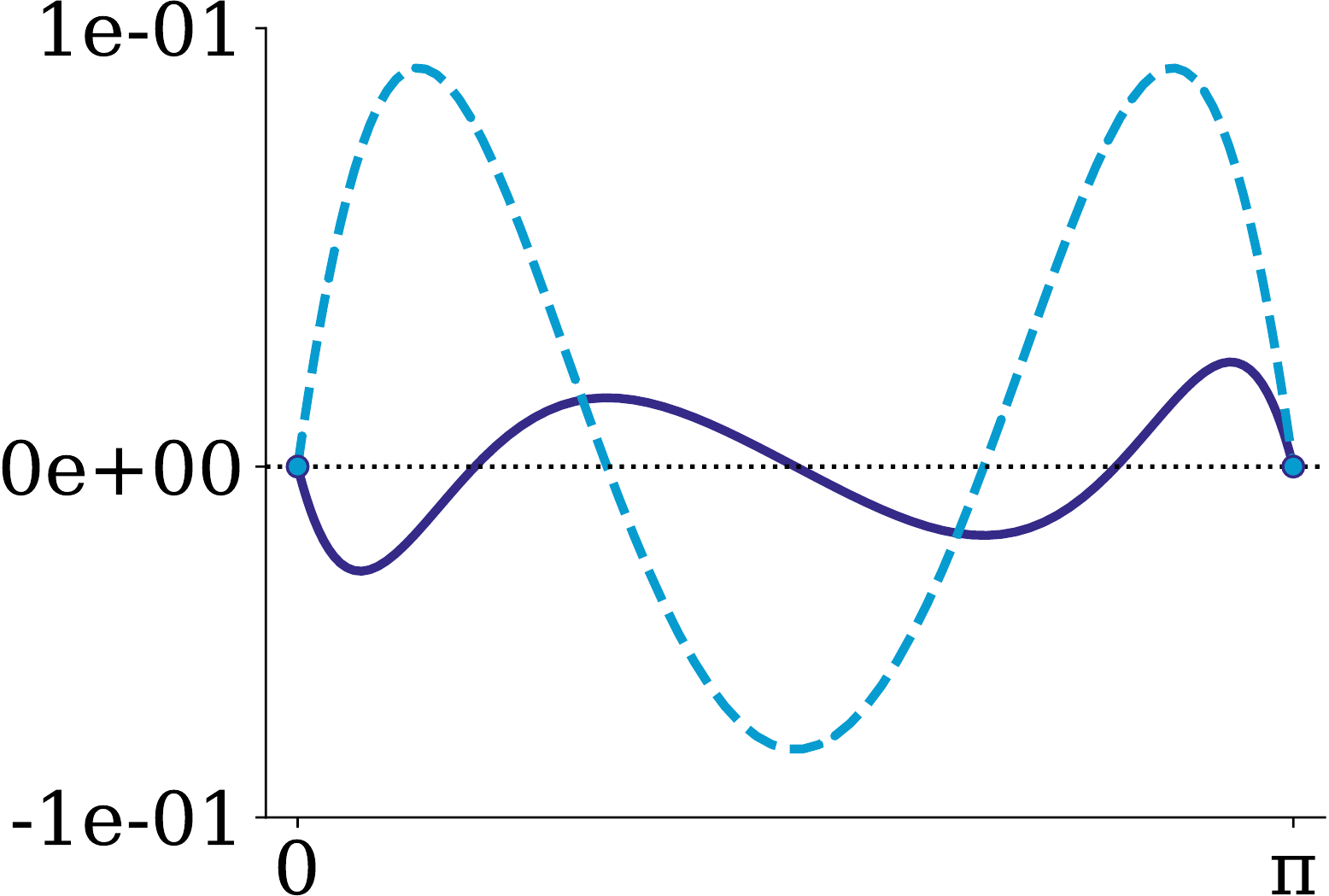} 
& 
\begin{overpic}[width=keepaspectratio,width=0.26\textwidth]{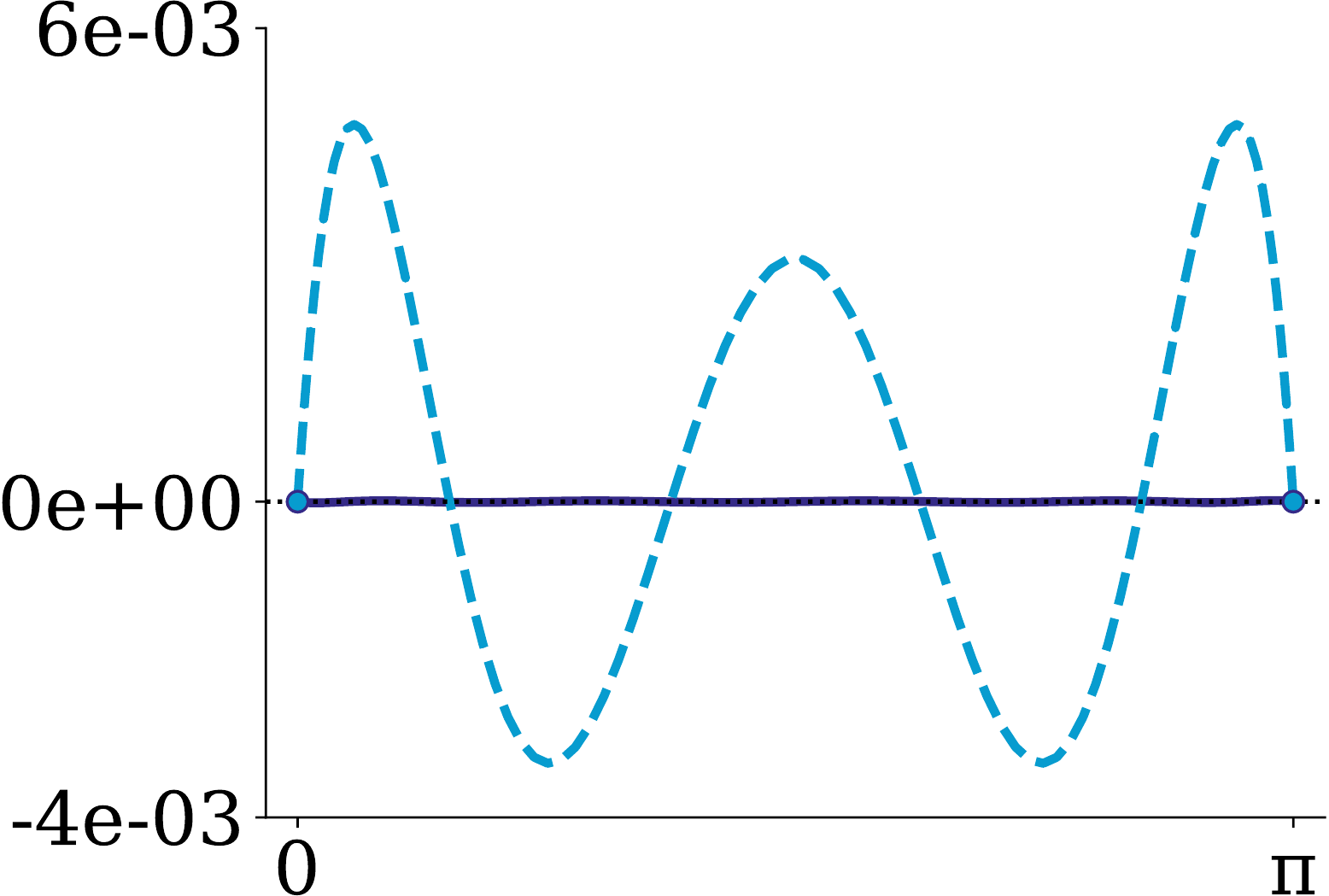} 
\end{overpic}
\\[10pt]
\hspace{40pt}

\begin{overpic}[
width=keepaspectratio,width=0.26\textwidth]{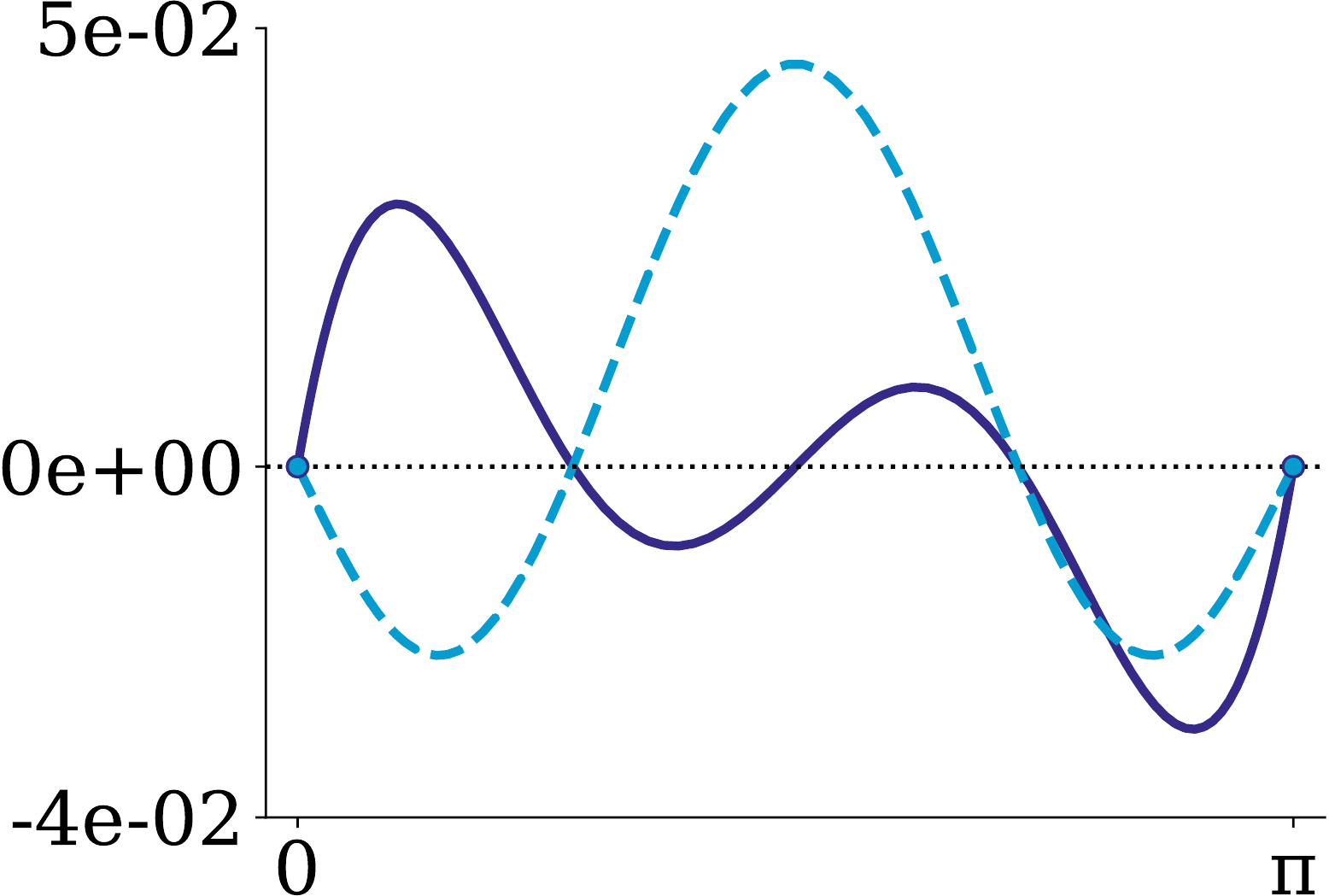} 
\put(40,-10){ \textbf{Initial} }
\put(-45,30){\fbox{$\boldsymbol{p = 3}$}}
\end{overpic} 
&
\begin{overpic}[width=keepaspectratio,width=0.26\textwidth]{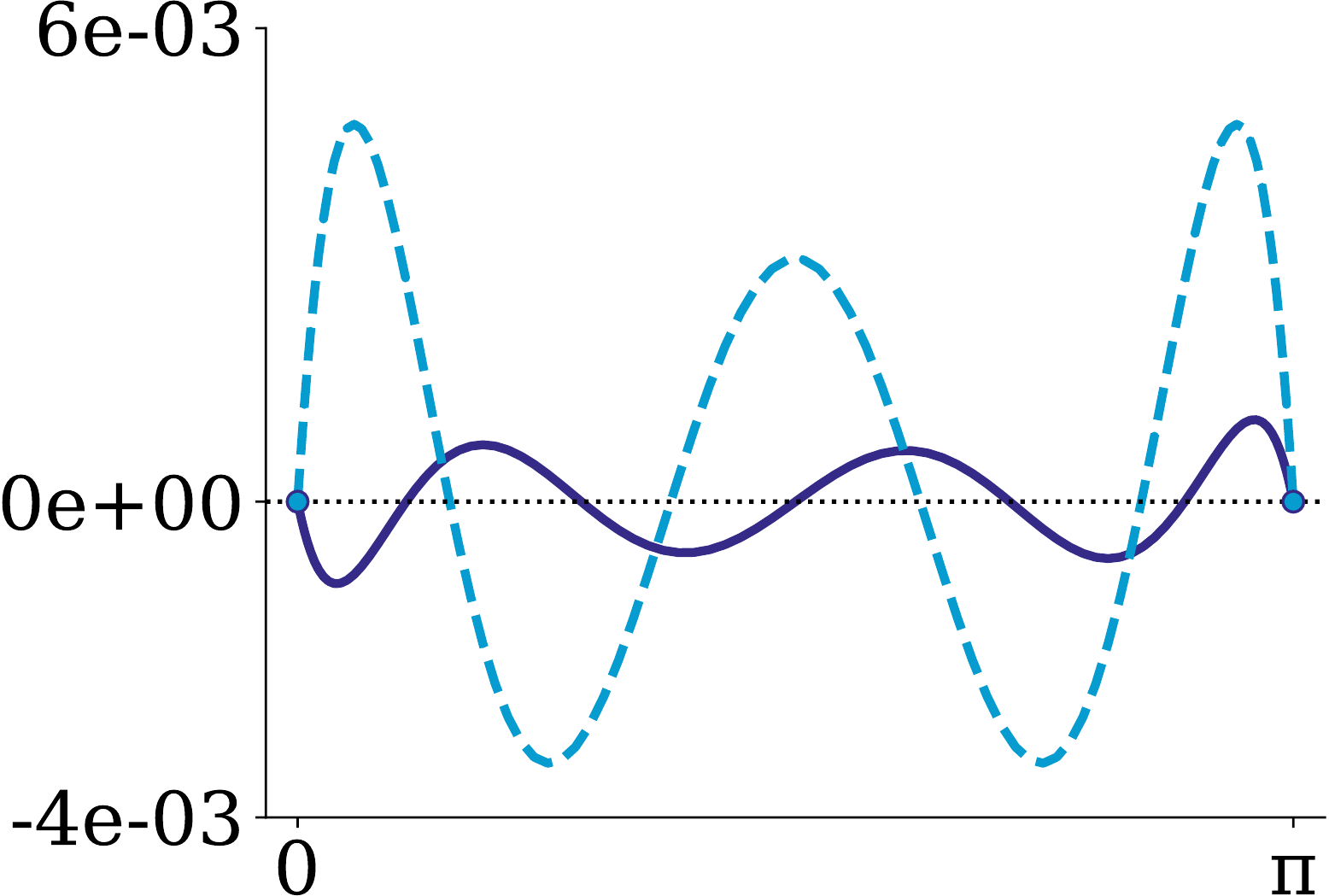} 
\put(35,-10){ $\boldsymbol{q = 2p-1}$ }
\end{overpic} 

& 
\begin{overpic}[width=keepaspectratio,width=0.26\textwidth]{plt/TN_OPT_Semi-Circle_P_3_Q_10_N_1} 
\put(35,-10){ $\boldsymbol{q = 10}$ }
\put(-20,-25){ \textbf{Optimal}}
\end{overpic}\\[20pt] 
\end{tabular}
\caption{\label{fig:tn_semicircle_1ele} tangent ($\boldsymbol t$) and normal ($\boldsymbol n$) error components approximating a semi-circle with one element before and after optimizing the internal nodes. 
}
\end{figure}
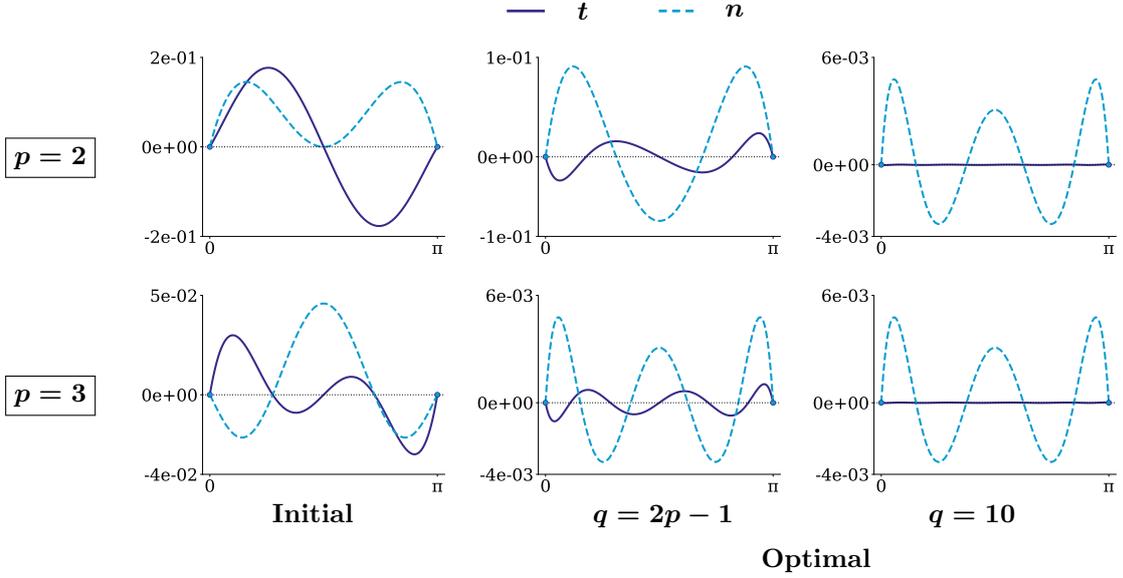

Let us now study the non-linear equations separately; during optimization, we impose zero tangent error (weakly) in $q-1$ equations. If we assume that solving the non-linear equations behaves like interpolation, we expect at least $q-1$ roots in the tangent error function plus the two roots at the interfaces. This is shown in Figure \ref{fig:tn_semicircle_1ele} for the $q=2p-1$ case: the optimized tangent error
has 5 and 7 roots for $p=2,3$, respectively.

We use the $2p-2$ remaining equations to impose the total error equal to zero. Assume that we can make the tangent error as small as desired by increasing $q$. Then, at the optimum, we can think of these $2p-2$ equations essentially imposing zero normal error (weakly). Following the reasoning we used for $s$, we can expect $2p$ roots along the normal component. Again, the extra two come from the interpolatory interfaces. This can be appreciated in the plots from the optimizad case in Figure \ref{fig:tn_semicircle_1ele}. We have $2p$ roots in the normal error and a tangent error that vanishes as we increase $q$.

\subsection{Discussion for 3D curves}

We have discussed the 2D case and how the optimal error behaves in terms of the tangent and normal components. We will now extend our results to the 3D case. In this case, the error,  $\be = \bx - \ba \circ s$, is decomposed as:
\begin{equation}\be = (\be\cdot \boldsymbol t)\boldsymbol t+(\be\cdot \boldsymbol n)\boldsymbol n +(\be\cdot \boldsymbol b)\boldsymbol b.
\end{equation}
Here, $\{\boldsymbol t,\boldsymbol n,\boldsymbol b\}$ are the curve tangent, normal and binormal vectors, respectively. 

Our physical mesh uses polynomials of degree $p$ in 3D space with fixed end-points. So, at each element, we have $3 (p-1)$ degrees of freedom. As for the 2D case, we impose in $q-1$ equations zero tangent error (weakly) and the $3(p-1)$ equations impose total zero error (weakly). Again, we assume that the tangent error decreases as $q$ increases. At the optimum, the combined solution implies that we have $3 (p-1)$ equations imposing zero along both the normal and binormal components. In analogy with the 2D discussion, we now expect at least $\lfloor\frac 32 (p-1)\rfloor$ interpolation points along each component: $\{\boldsymbol n, \boldsymbol b\}$. Adding the two end-points gives $(\lfloor\frac 32(p-1)\rfloor + 2)$ roots per  component. 

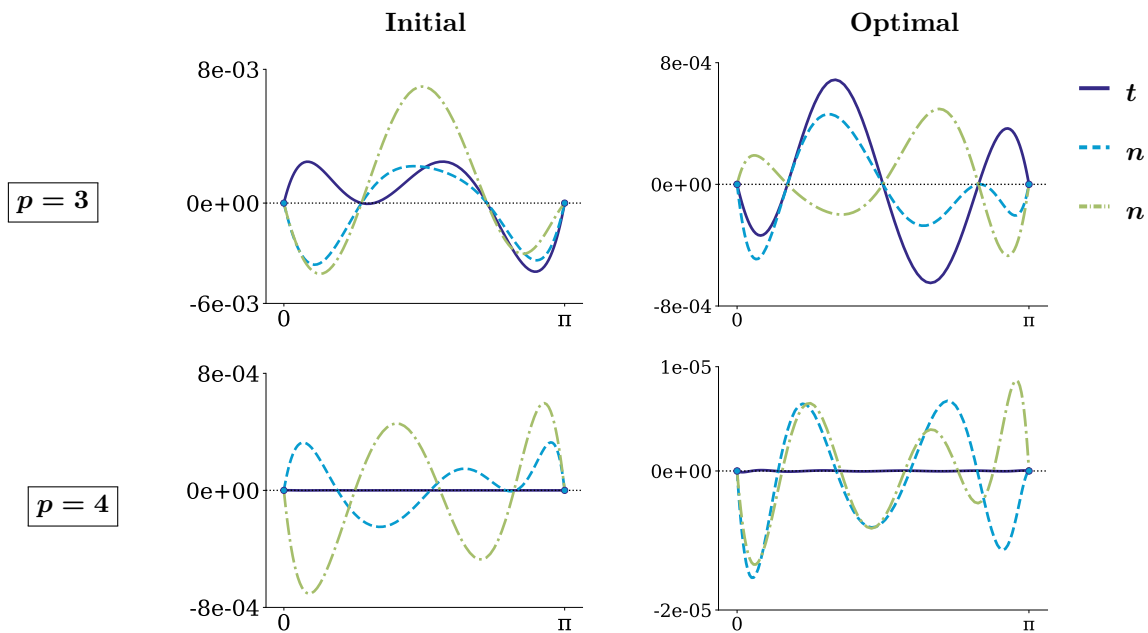
\begin{figure}[h!]
\vspace{20pt}
\flushleft
\begin{tabular}{cc} 
\hspace{50pt}
\begin{overpic}[width=keepaspectratio,width=0.34\textwidth]{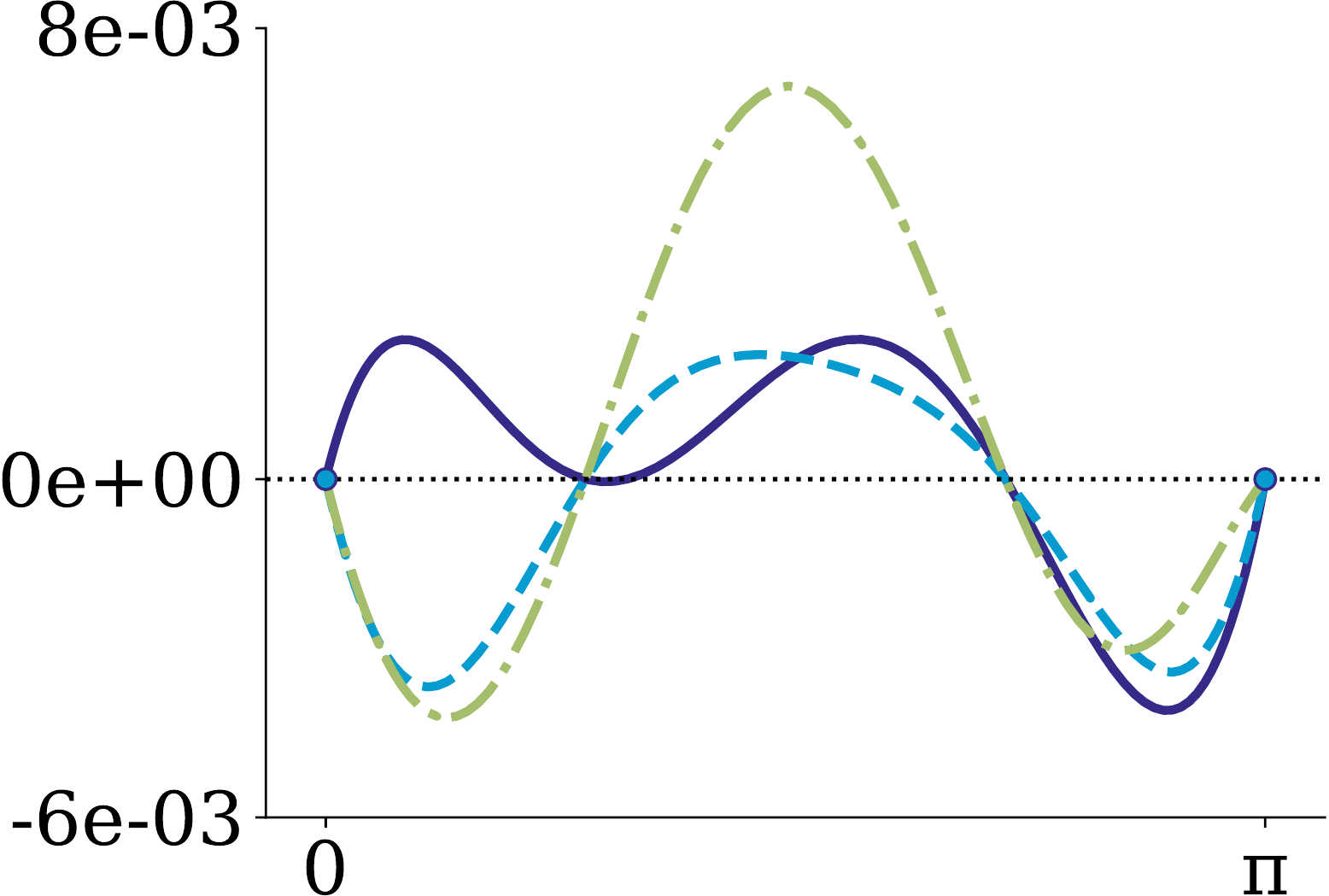}  
\put(-45,30){\fbox{$\boldsymbol{p = 3}$}}
\put(50,75){\textbf{Initial}}
\end{overpic}
&
\hspace{10pt}
\begin{overpic}
[width=keepaspectratio,width=0.34\textwidth]{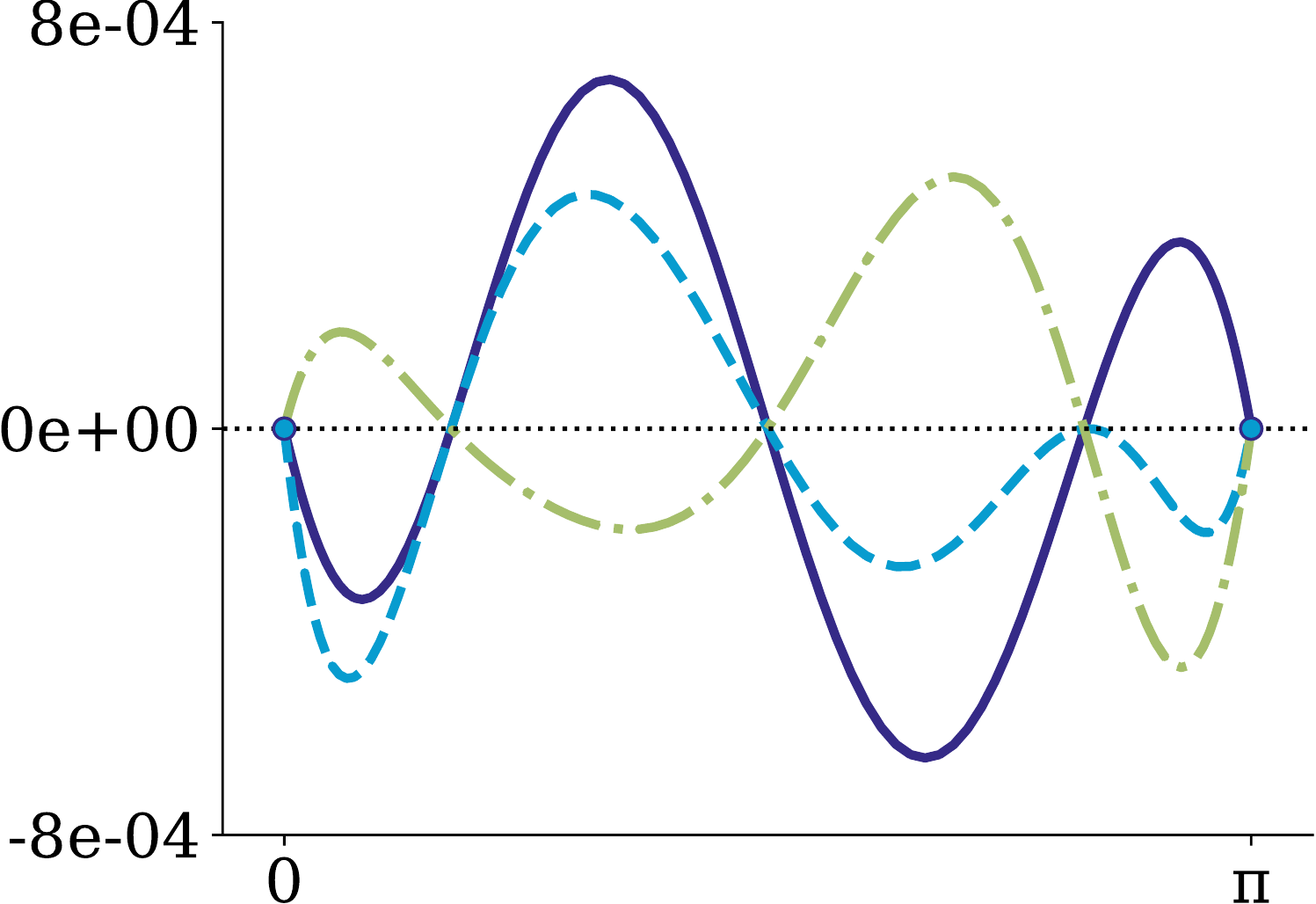} 
\put(50,75){\textbf{Optimal}}
\put(105,60){ \tikz{\draw (0.0,0.0)--(0.4,0.0) [par_1, line width=1.5pt] }}
\put(115,55){ \tikz{\draw (0.8,0.0)  node {$\boldsymbol t$} }}

\put(105,45){ \tikz{\draw (0.0,0.0)--(0.4,0.0) [densely dashed, par_4,line width=1.5pt] }}
\put(115,40){ \tikz{\draw (0.8,0.0)  node {$\boldsymbol n$} }}

\put(105,30){ \tikz{\draw (0.0,0.0)--(0.4,0.0) [densely dashdotted, par_7, line width=1.5pt] }}
\put(115,25){ \tikz{\draw (0.8,0.0)  node {$\boldsymbol n$} }}
\end{overpic}
\\[10pt]
\hspace{50pt}
\begin{overpic}[width=keepaspectratio,width=0.34\textwidth]{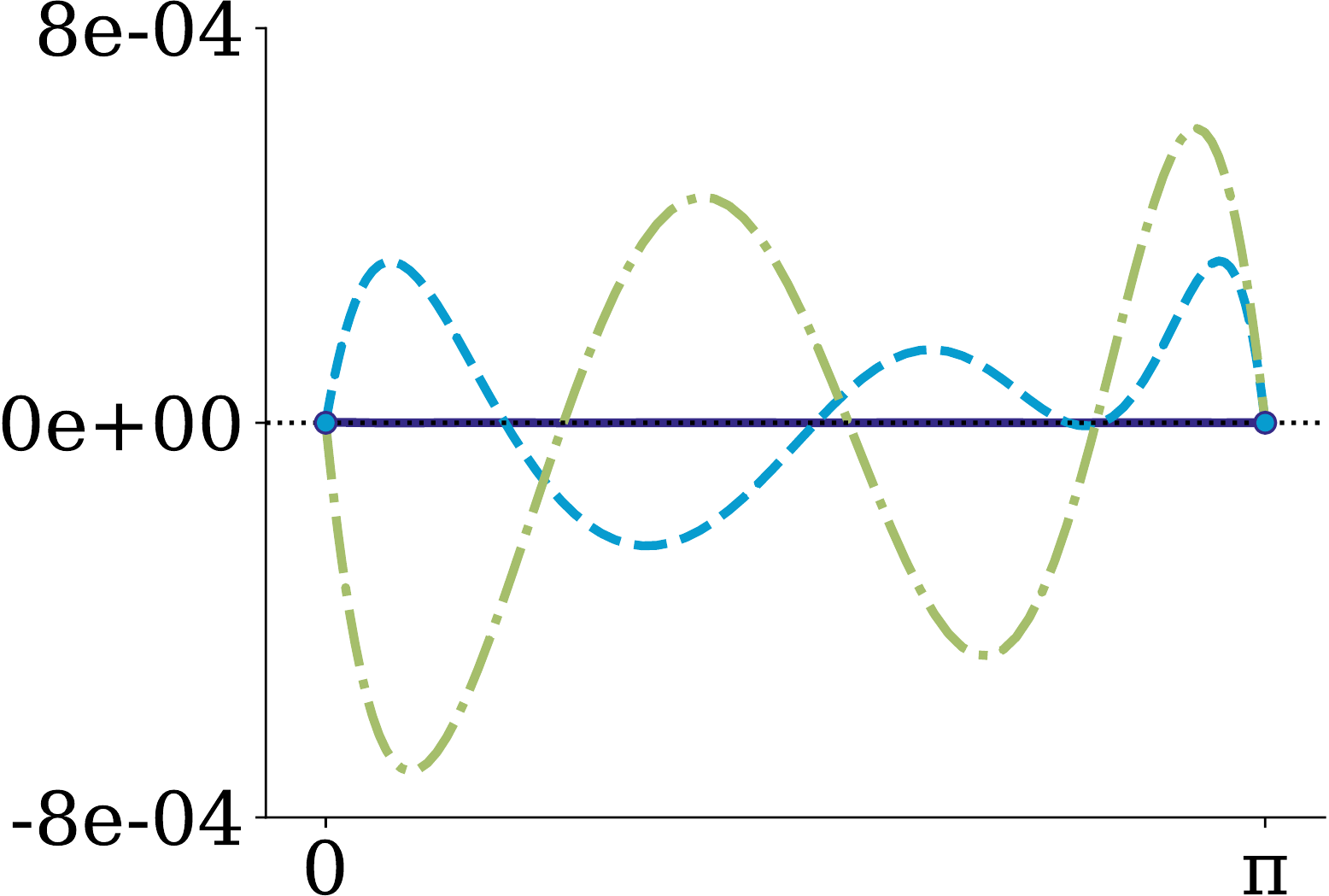} 
\put(-40,30){\fbox{$\boldsymbol{p = 4}$}}
\end{overpic}& 
\hspace{10pt}
\begin{overpic}[
width=keepaspectratio,width=0.34\textwidth]{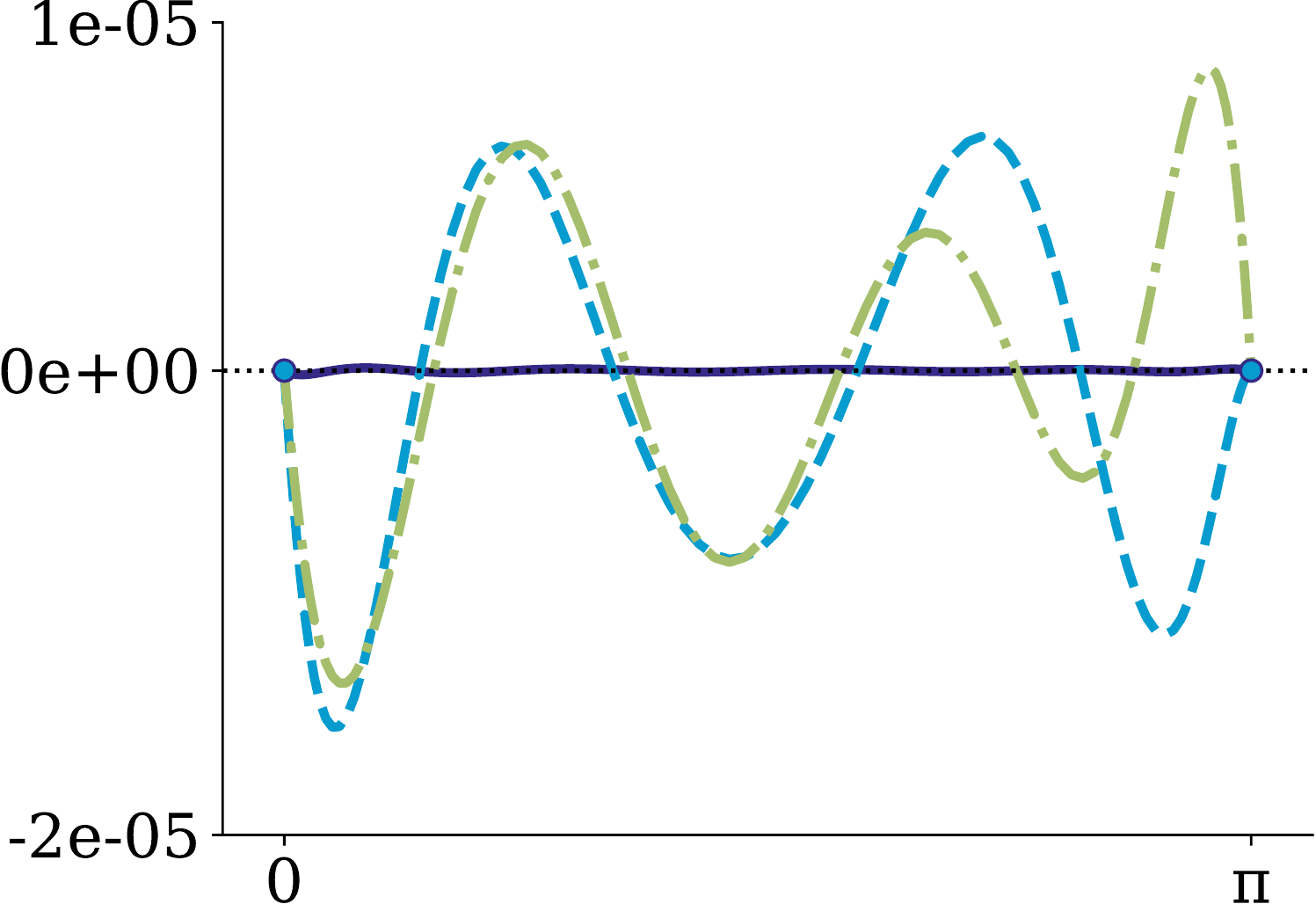} 
\end{overpic}
\end{tabular}
\caption{\label{fig:tnb_sphere_1ele} tangent,  normal and binormal $\{\boldsymbol t,\boldsymbol n,\boldsymbol b\}$  point-wise errors  approximating a sphere arc with one element before and after minimizing the constrained disparity.}
\end{figure}

In Figure \ref{fig:tnb_sphere_1ele}, we show the error plots before and after optimizing the constrained disparity approximating a sphere arc with a single element. As for the 2D case, as $\tilde \bx$ follows the image of $\ba\circ \tilde s$,  $\tilde s$ minimizes the tangent error. Notice that when $s$ is of degree $q=10$, the tangent error is negligible compared to the other two components. Also, observe how we obtain both along the normal and binormal directions: $5 = \lfloor\frac 32(3-1)\rfloor + 2$ roots for $p=3$ and at least $6 = \lfloor\frac 32(4-1)\rfloor + 2$ for $p=4$.

\section{Numerical experiments on CAD data} \label{sec:numerical_results}
Here we study the performance of the new solver both in terms of accuracy and computational times. We use CAD models from the ESP database, evaluate the curve and derivatives using the EGADS \emph{Julia} wrap, and run the HPC experiments on the \emph{BSC MareNostrum 4} supercomputer. 

\begin{figure}[h]
\flushleft
\vspace{10pt}
\hspace{20pt}
 \begin{overpic}[width = keepaspectratio, width = 0.2\textwidth] {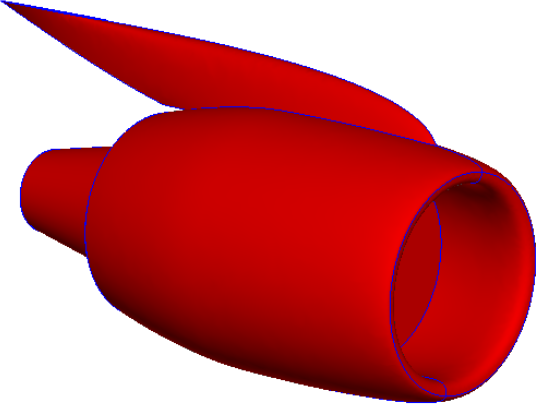} 
 \put(110,30){
 \scalebox{0.9}{\begin{tabular}{ccc}
  & \textbf{Unconstrained} & \textbf{Constrained}\\
Iterations & 823&202\\
Line Searches& 2296&  640\\
$|\nabla E|<1.e^{-12}$&4&29\\
CPU minutes&14.5&2.5
 \end{tabular}}}
 \end{overpic}

 \vspace{30pt}
 \begin{tabular}{ccc}
 \mc{3}{\hspace{45pt}
\begin{overpic}[width = keepaspectratio, width = 0.6\textwidth]{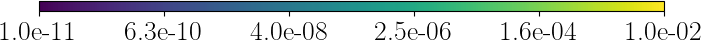}
\put(-20,3){\textbf{$\boldsymbol{|error|}$ }}  
  \end{overpic}   }\\[10pt]  
  
\begin{overpic}[width=keepaspectratio,  width = 0.3\textwidth]
    {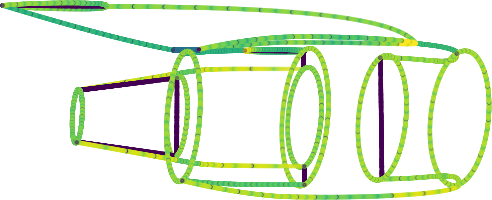}
    \put(50,-15){\textbf{Initial}}
    \end{overpic}
   &
    \begin{overpic}[width=keepaspectratio, width = 0.3\textwidth]
        {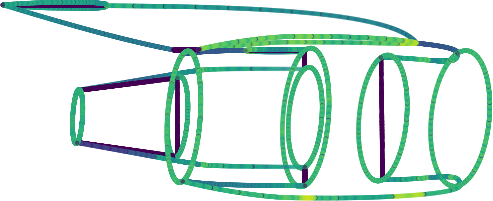}
        \put(20,-15){\textbf{Unconstrained}}
    \end{overpic} 
    & 
    \begin{overpic}[width=keepaspectratio, width = 0.3\textwidth]{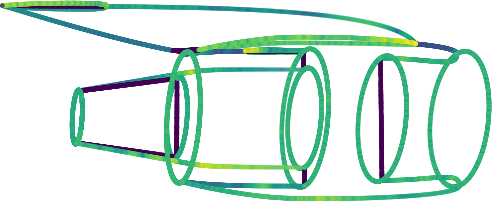}
    \put(30,-15){\textbf{Constrained}}
    \end{overpic}\\[15pt]

  \end{tabular}
\caption{\label{fig:nacelle_fix_free} optimization results from approximating the 37 edges of an aircraft nacelle using meshes with 12 elements (per curve) and $p=2$. Top: number of iterations, line searches, curves that converged and CPU times. Bottom: the point-wise errors.}
\end{figure}

In Figure \ref{fig:nacelle_fix_free} we show the results from approximating the curves of a nacelle body using unconstrained and constrained optimizations. Both solutions reduce the approximation errors, with a greater reduction in those curves that are planar. The constrained disparity requires $\sim 75 \%$ fewer iterations and reduces the computing time by a factor of 5.7. Moreover, the underlying CAD parametrization consists of B-splines with potential non-smooth derivatives, affecting the optimizer, especially for the unconstrained version: the number of line searches is very high, affecting the convergence. Only four curves were able to reach the stop criteria: $|\nabla E|<1.e^{-12}$. On the other hand, the constrained problem converged for 29 out of 37 curves.

\subsection{Parallel element optimization}
\begin{figure}[h!]
\centering
 \begin{tabular}{cc}
 
  \includegraphics[width=keepaspectratio, width = 0.45\textwidth]{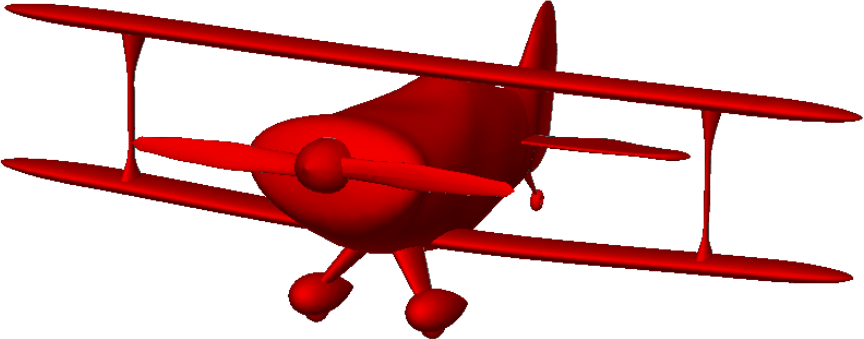} & 
  \includegraphics[width=keepaspectratio, width = 0.45\textwidth]{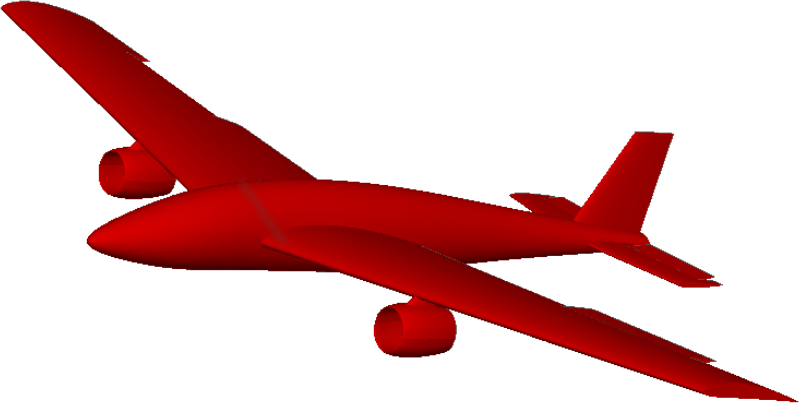} \\
  \includegraphics[width=keepaspectratio, width = 0.45\textwidth]{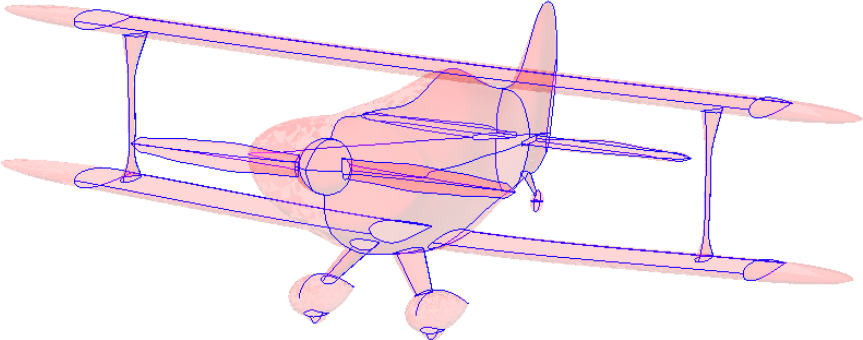} &
  \includegraphics[width=keepaspectratio, width = 0.45\textwidth]{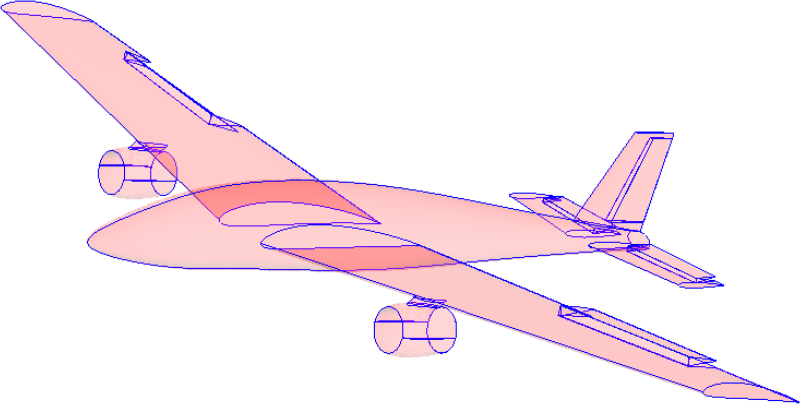} \\
 \textbf{Aircraft 1: 141 curves} &\textbf{Plane 2: 86 curves} 
 \end{tabular} 
  \caption{\label{fig:plane_models} two aircraft CAD prototypes (top) highlighting all the curves (bottom).}
 \end{figure}

\begin{figure}[h!]
\centering
  \begin{tabular}{cc}
\hspace{10pt}
 \begin{overpic}[width=keepaspectratio, width = 0.4\textwidth]{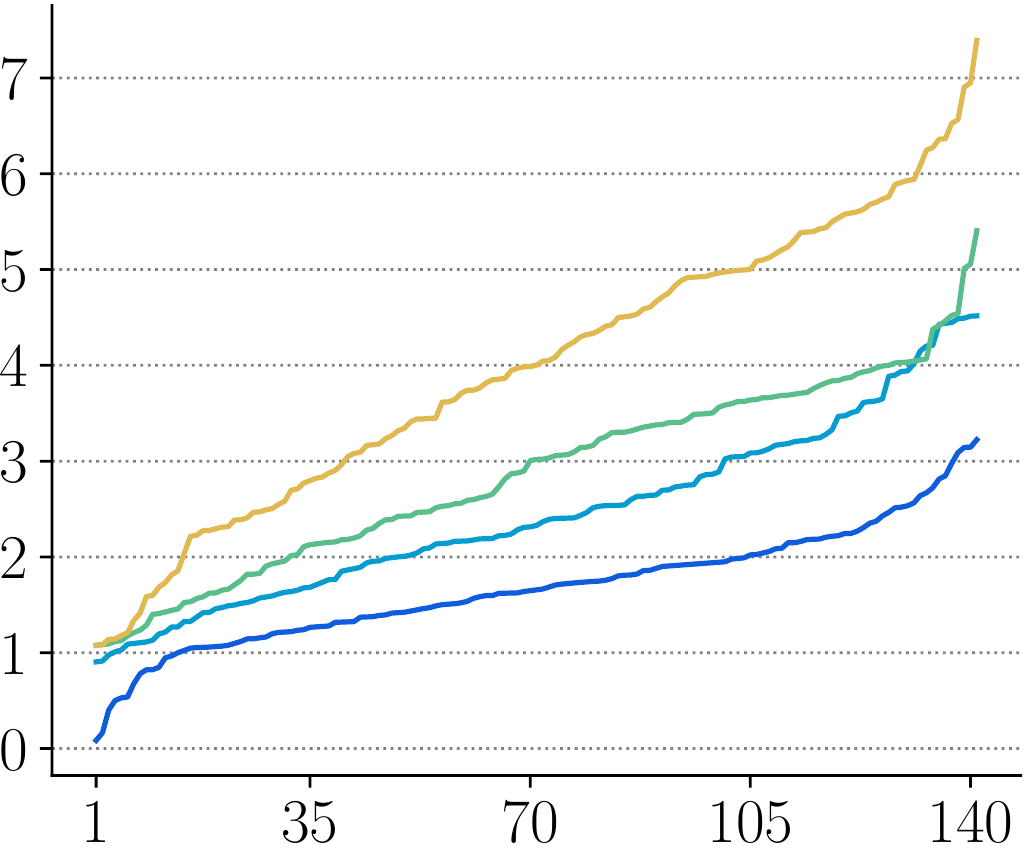}
 \put(15,90){\textbf{Aircraft 1, $\boldsymbol{p=2}$}}
 \put(-15,20){\rotatebox{90}{\textbf{speedup:} $\boldsymbol{\frac{t_s}{t_p}}$}}
 \put(70,79){\textcolor{par_1}{\scriptsize{\textbf{CPUs:}}}}
  \put(95,79){\textcolor{par_1}{$\scriptstyle{\boldsymbol{48}}$}}
  \put(95,62){\textcolor{par_1}{$\scriptstyle{\boldsymbol{36}}$}}
  \put(95,53){\textcolor{par_1}{$\scriptstyle{\boldsymbol{24}}$}}
  \put(95,41){\textcolor{par_1}{$\scriptstyle{\boldsymbol{12}}$}}
  \put(10,60){\textcolor{par_1}{$\scriptstyle{\boldsymbol{average \left( \frac{t_s}{t_{48} } \right)\ =\ 4}}$} }
 \end{overpic}
&
\hspace{20pt}
 \begin{overpic}[width=keepaspectratio, width = 0.4\textwidth]{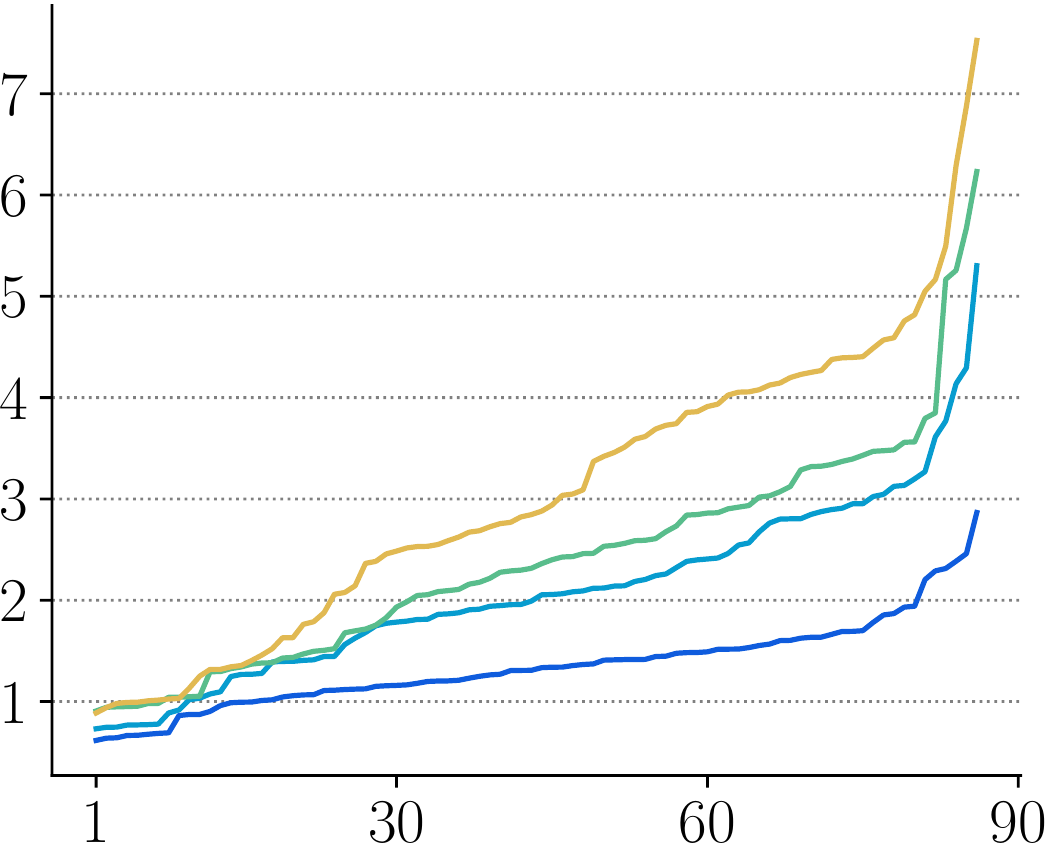}
 \put(15,90){\textbf{Aircraft 2, $\boldsymbol{p=3}$}}
 \put(70,79){\textcolor{par_1}{\scriptsize{\textbf{CPUs:}}}}
  \put(94,77){\textcolor{par_1}{$\scriptstyle{\boldsymbol{48}}$}}
  \put(94,64){\textcolor{par_1}{$\scriptstyle{\boldsymbol{36}}$}}
  \put(94,54){\textcolor{par_1}{$\scriptstyle{\boldsymbol{24}}$}}
  \put(94,32){\textcolor{par_1}{$\scriptstyle{\boldsymbol{12}}$}}
  \put(10,57){\textcolor{par_1}{$\scriptstyle{\boldsymbol{average \left( \frac{t_s}{t_{48} } \right)\ =\ 3}}$} }
 \end{overpic}\\[20pt]
  
  \hspace{10pt}
\begin{overpic}[width=keepaspectratio, width = 0.4\textwidth]{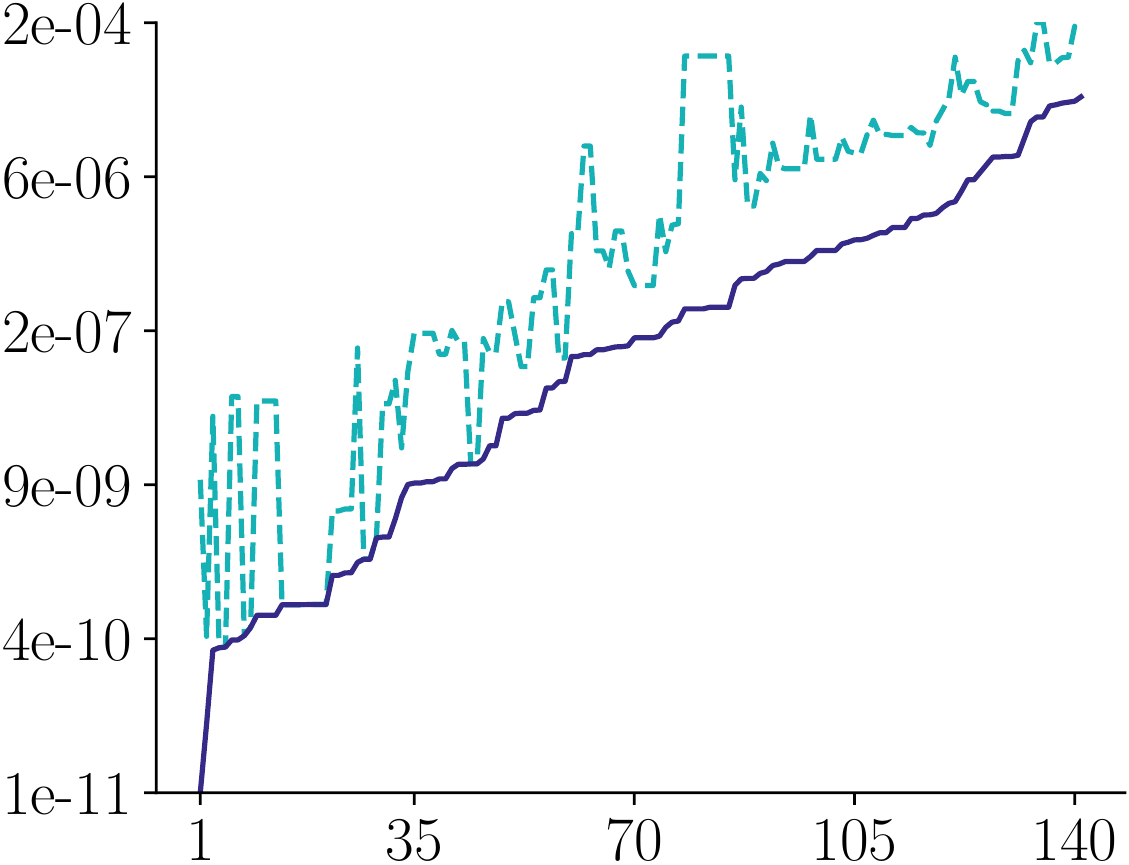}
 \put(-13,10){\rotatebox{90}{\textbf{disparity: }$\boldsymbol{||\cdot||_\sigma}$}}
 \put(65,20){$\scriptstyle{||\bx-\ba||_\sigma}$}
 \put(50,20){
 \tikz{
 \draw(0.0,0.0)--(0.4,0.0)[par_5, densely dashed, line width = 1.0pt]}}
 \put(65,10){$\scriptstyle{||\bx^\star-\ba\circ s^\star||_\sigma}$}
 \put(50,10){
 \tikz{
 \draw(0.0,0.0)--(0.4,0.0)[par_1,  line width = 1.0pt]}}
 \put(40,-10){\textbf{Curves}}
\end{overpic}
&
  \hspace{20pt}
\begin{overpic}[width=keepaspectratio, width = 0.4\textwidth]{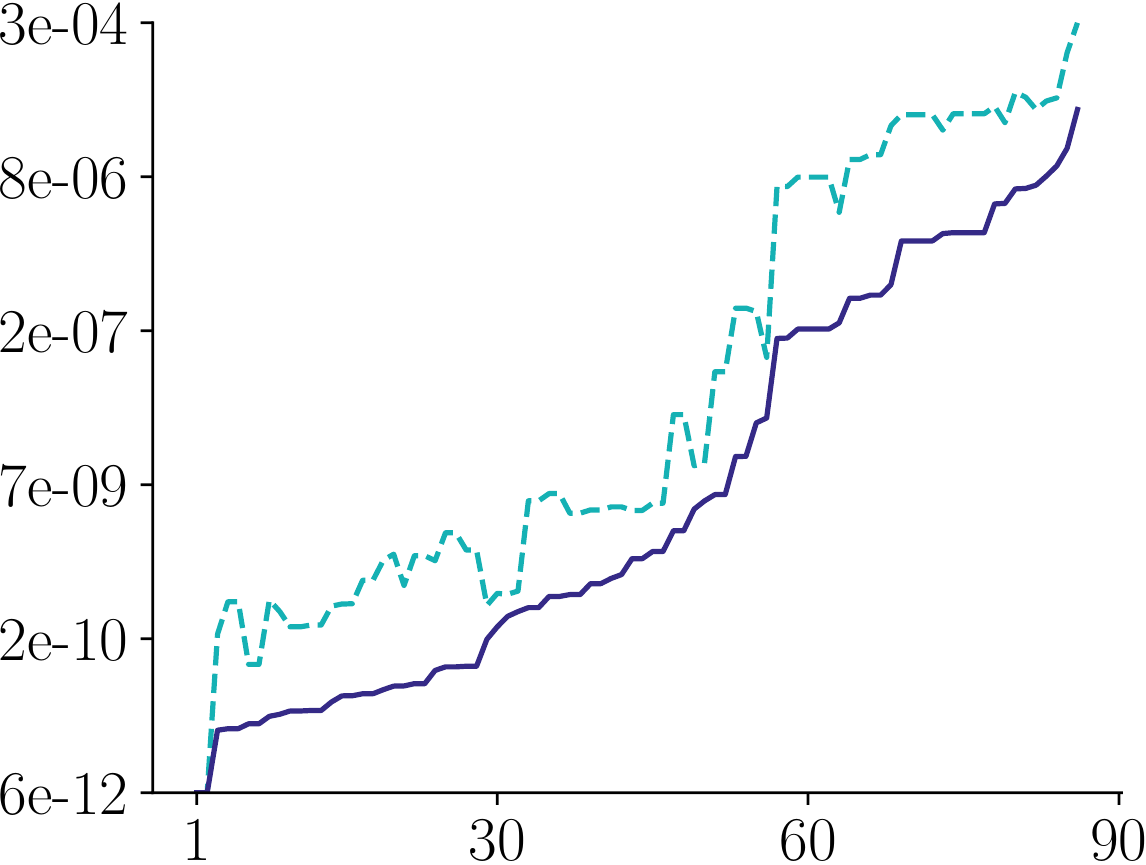}
 \put(65,20){$\scriptstyle{||\bx-\ba||_\sigma}$}
 \put(50,20){
 \tikz{
 \draw(0.0,0.0)--(0.4,0.0)[par_5, densely dashed, line width = 1.0pt]}}
 \put(65,10){$\scriptstyle{||\bx^\star-\ba\circ s^\star||_\sigma}$}
 \put(50,10){
 \tikz{
 \draw(0.0,0.0)--(0.4,0.0)[par_1,  line width = 1.0pt]}}
 \put(40,-10){\textbf{Curves}}
\end{overpic}\\[20pt]
  
 \end{tabular} 
  \caption{\label{fig:element_hpc} parallel optimization results for the two aircrafts showing the speedup factors at each curve (top), and the initial and final constrained disparities (bottom).}
 \end{figure}

Now we carry out a parallel study by fixing the number of elements and increasing the number of cores. We use the two aircraft prototypes from Figure \ref{fig:plane_models}. In Figure \ref{fig:element_hpc} we show the CPU performance after optimizing meshes made of 48 elements (per curve) and varying polymial degree. Regarding accuracy (bottom plots), we see that the optimized meshes improve the error in the approximation. On average, the error relative to the original approximation improves by 81\% for $p=2$ and 85 \% for the $p=3$ case. 

The CPU speedup factors are shown in Figure \ref{fig:element_hpc} (top). Increasing the number of cores does not lead to a significant computational gain. Using 48 cores, we obtain, on average, a speedup of four in $p=2$ and three in $p=3$. Also, many curves have a value of one or less, meaning slower computing times. It can be explained by looking at the workload: while most cores have converged, the few that remain computing increase the overall times. 

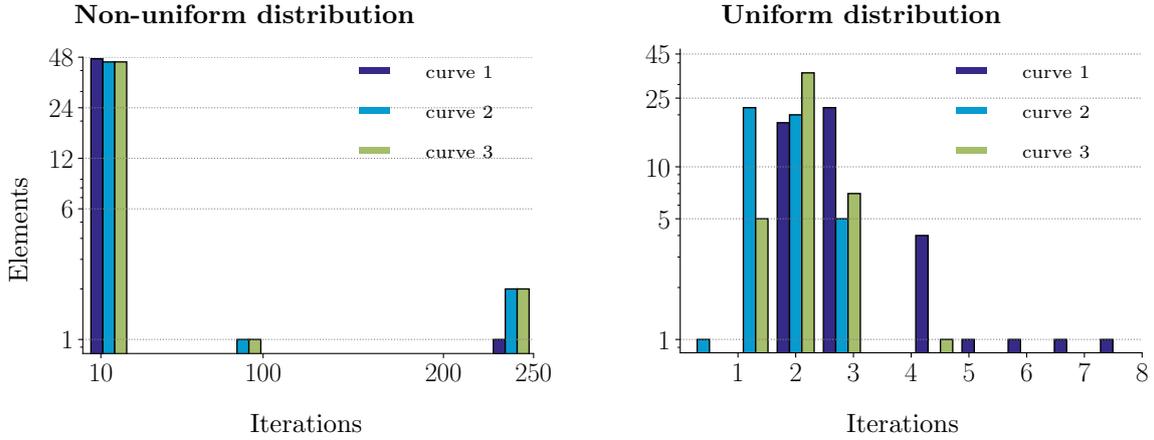
\begin{figure}[h!]
 \flushright
 \vspace{26pt}
 \begin{tabular}{cc}
   \begin{overpic}[width = keepaspectratio, width = 0.43 \textwidth]{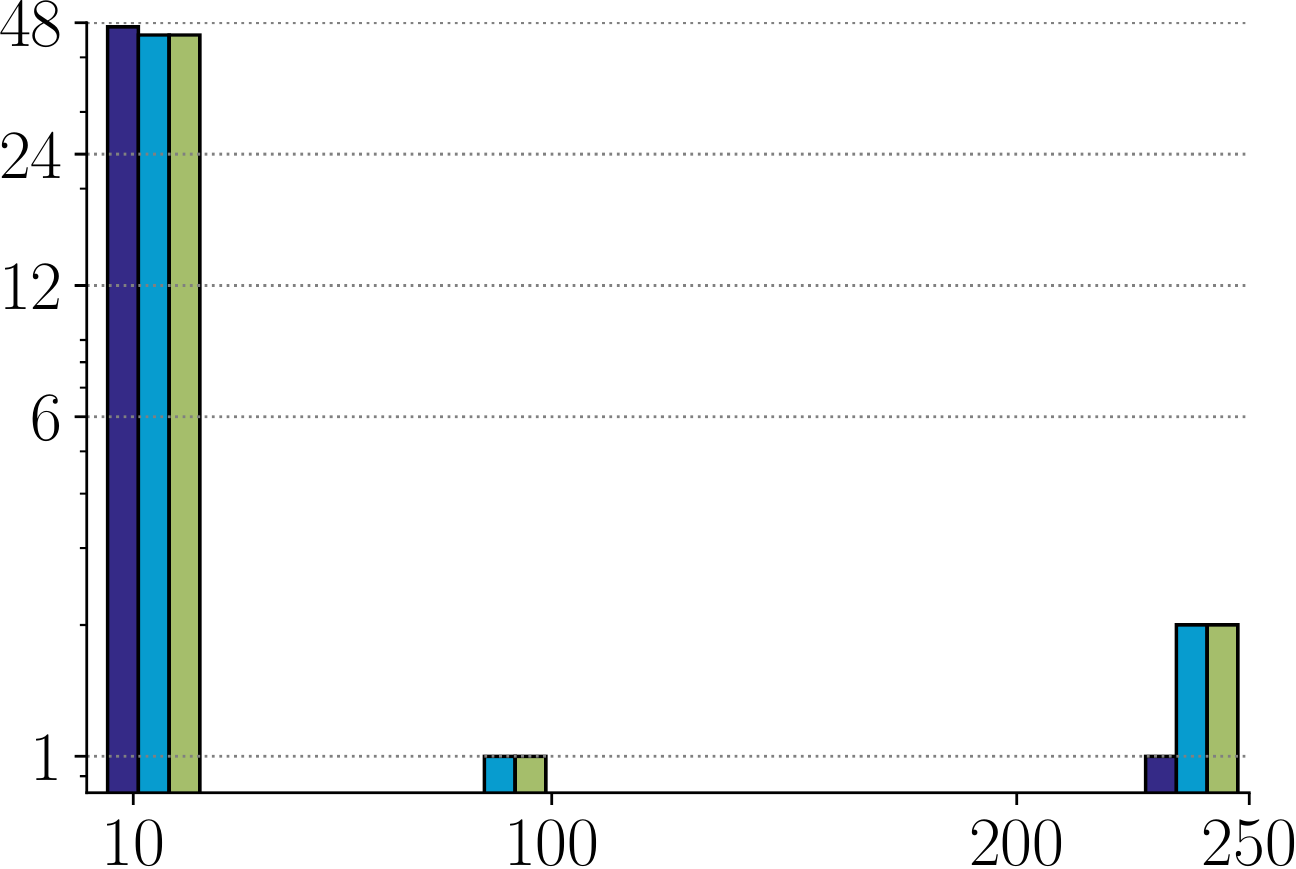}
  \put(5,72){\textbf{Non-uniform distribution} }
  \put(-8,20){\rotatebox{90}{Elements} } 
    \put(40,-10){Iterations} 
  
  \put(60,62) {
  \tikz{  \fill (0.0,0.0) -- (0.4,0.0) --(0.4,0.08)--(0.0,0.08)--cycle[par_1]}  }
  \put(75,61){\scriptsize{curve 1}}
  
  \put(60,54) {
  \tikz{  \fill (0.0,0.0) -- (0.4,0.0) --(0.4,0.08)--(0.0,0.08)--cycle[par_4]}  }
  \put(75,53){\scriptsize{curve 2}}
  
  \put(60,46) {
  \tikz{  \fill (0.0,0.0) -- (0.4,0.0) --(0.4,0.08)--(0.0,0.08)--cycle[par_7]}  }
   \put(75,45){\scriptsize{curve 3}}
   
     \end{overpic} & \hspace{20pt}
  \begin{overpic}[width = keepaspectratio, width = 0.43 \textwidth]{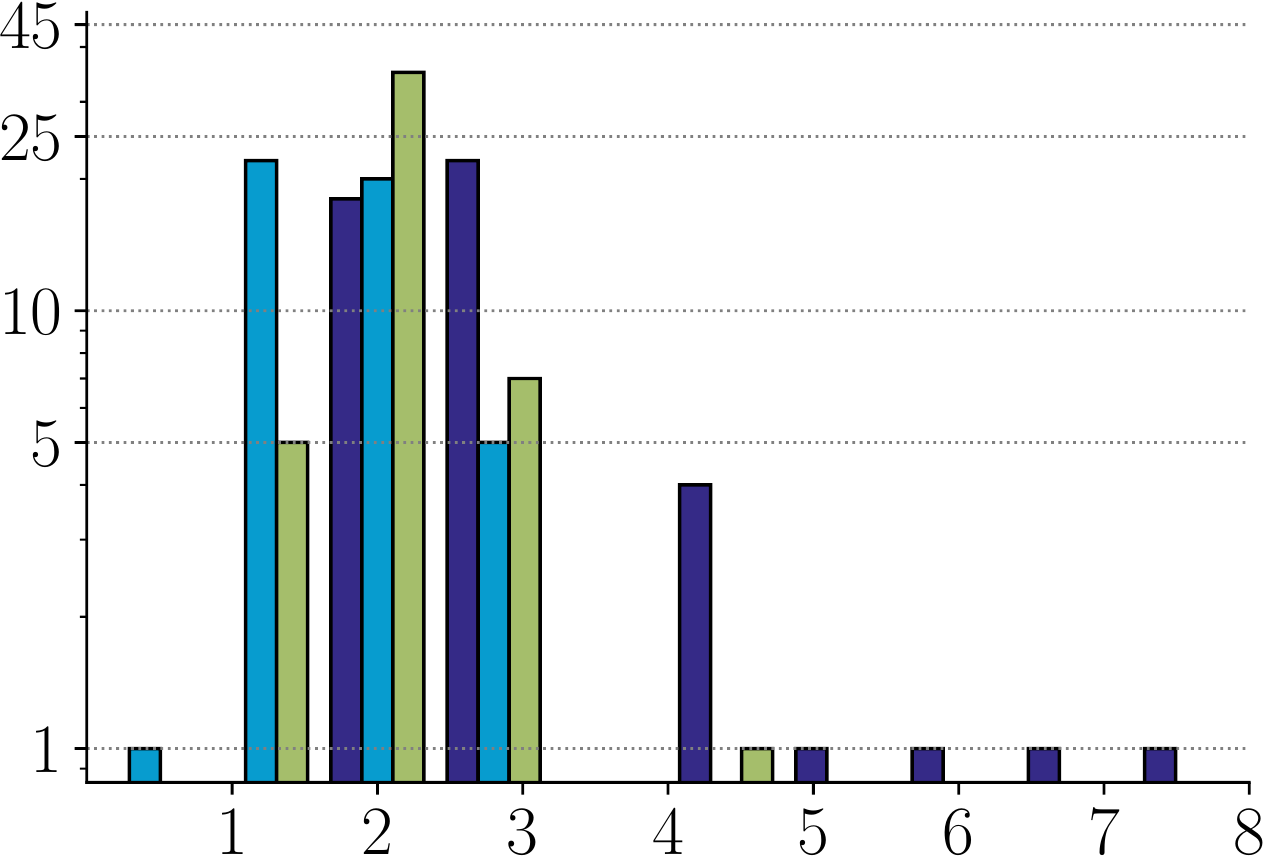}
  \put(15,72){\textbf{Uniform distribution} }
  
  \put(60,62) {
  \tikz{  \fill (0.0,0.0) -- (0.4,0.0) --(0.4,0.08)--(0.0,0.08)--cycle[par_1]}  }
  \put(75,61){\scriptsize{curve 1}}
  
  \put(60,54) {
  \tikz{  \fill (0.0,0.0) -- (0.4,0.0) --(0.4,0.08)--(0.0,0.08)--cycle[par_4]}  }
  \put(75,53){\scriptsize{curve 2}}
  
  \put(60,46) {
  \tikz{  \fill (0.0,0.0) -- (0.4,0.0) --(0.4,0.08)--(0.0,0.08)--cycle[par_7]}  }
   \put(75,45){\scriptsize{curve 3}}
  
    \put(40,-10){Iterations} 
     \end{overpic}
     \end{tabular}\\[10pt]
     \caption{\label{fig:element_iters} distribution of the iterations per element taken by the optimizer over curves on the aircraft 1 model split between very large differences (left) vs. small differences (right). }
 \end{figure}

In Figure \ref{fig:element_iters} we show a histogram of the iterations at six curves from aircraft 1 separated by non-uniform and uniformly distributed iterations. For the left plot (non-uniform), we see that at least 46 out of 48 elements converged in less than 10 iterations but a couple of elements took more than 100. On the right, we see uniform iterations spread across elements, all converging in less than 10. 
 
\subsection{Parallel curve optimization}
\begin{figure}[h!]
 \centering
 \begin{tabular}{cc}
\raisebox{-55pt}{ \hspace{5pt} \begin{overpic}[width = keepaspectratio, width = 0.4\textwidth]{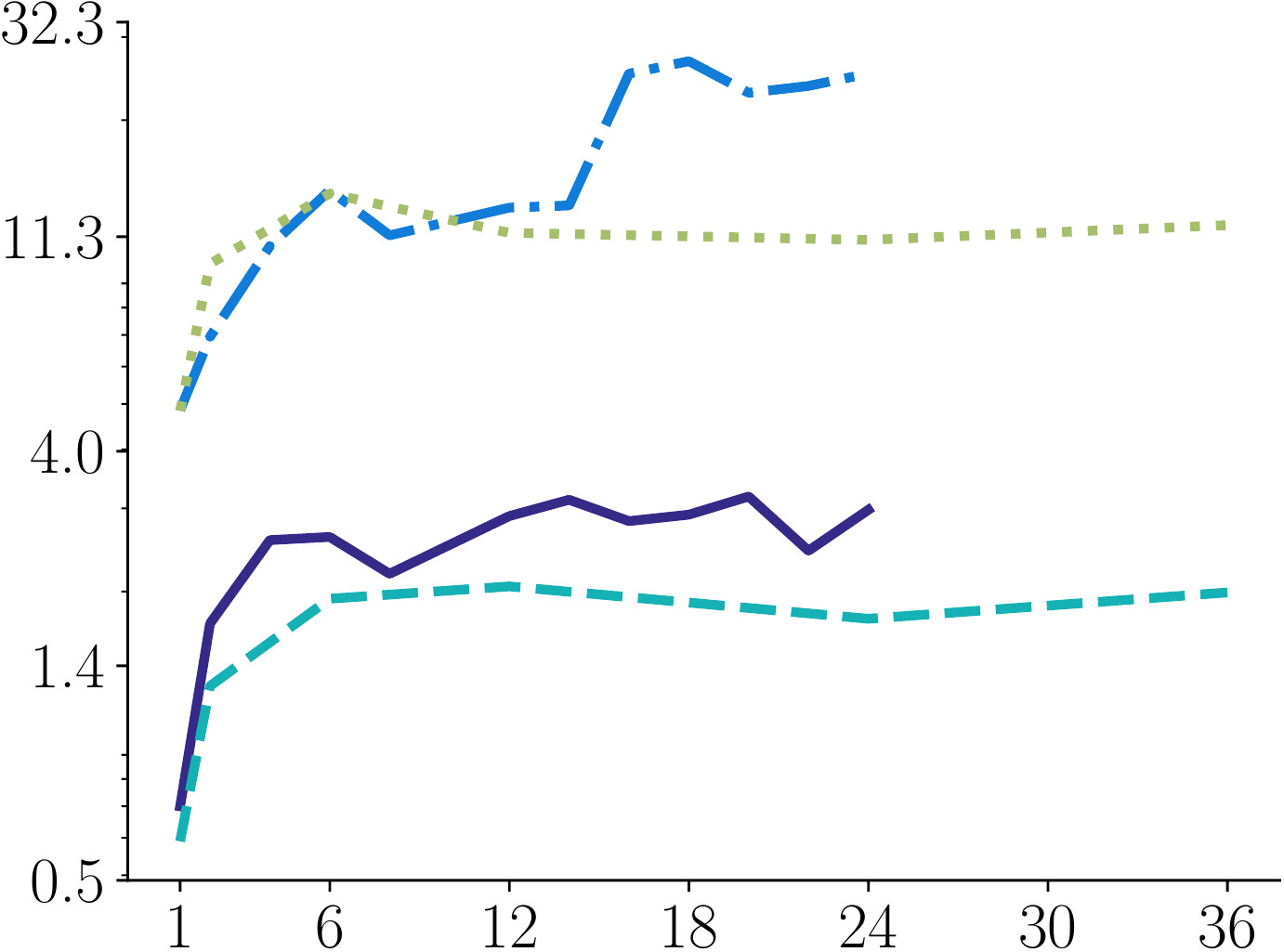}
  \put(-8,18){\scriptsize{\rotatebox{90}{\textbf{CPU time (s)}}}}
  \put(45,-7){\scriptsize{\textbf{Cores}}}
  \put(60,70){\scriptsize{$p = 3,$ 24 curves}}
  \put(65,50){\scriptsize{$p = 3,$ 1 curve}}
  \put(60,37){\scriptsize{$p = 2,$ 24 curves}}
  \put(65,20){\scriptsize{$p = 2,$ 1 curve}}
  \end{overpic}}
& \hspace{35pt}
 \begin{tabular}{c}
  \begin{overpic}[width = keepaspectratio, width = 0.35\textwidth]{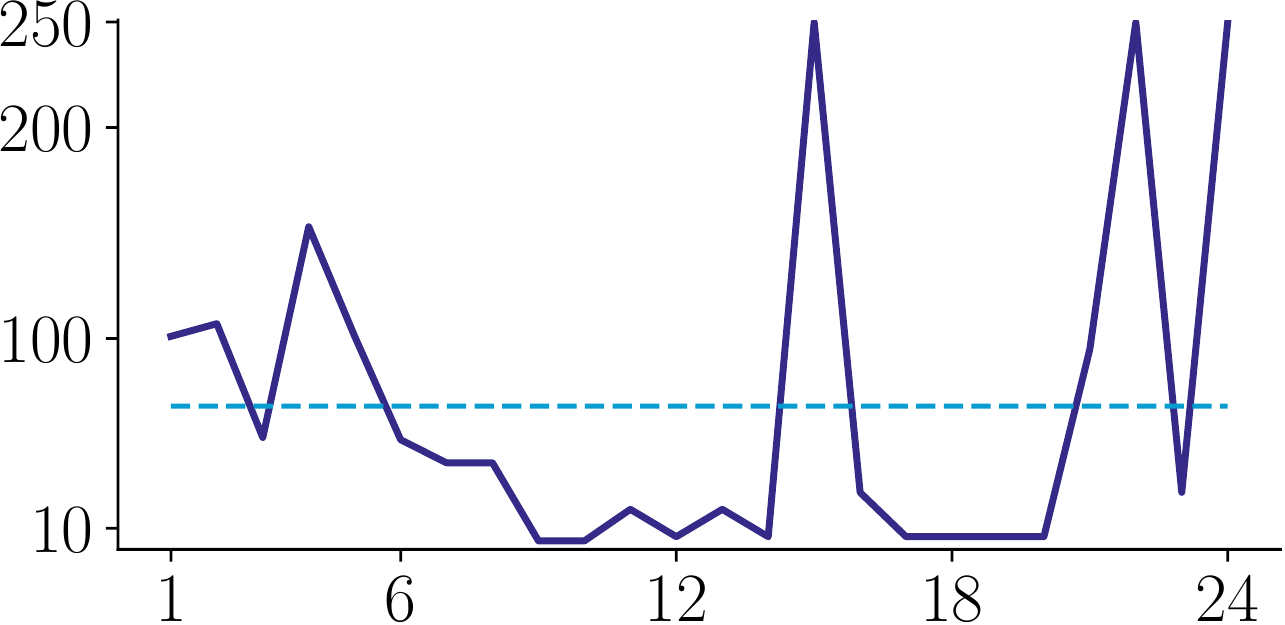}
  
  \put(-8,10){\scriptsize{\rotatebox{90}{\textbf{Iterations}}}}
  \put(36,18){
  \tikz{  \draw[draw=white,fill=white] (0,0) rectangle (2,0.4) }}
  \put(40,20){\scriptsize{\textcolor{par_4}{\textbf{average = 68}}}}
  \end{overpic}\\
  \hspace{3pt}
\begin{overpic}[width = keepaspectratio, width = 0.35\textwidth]{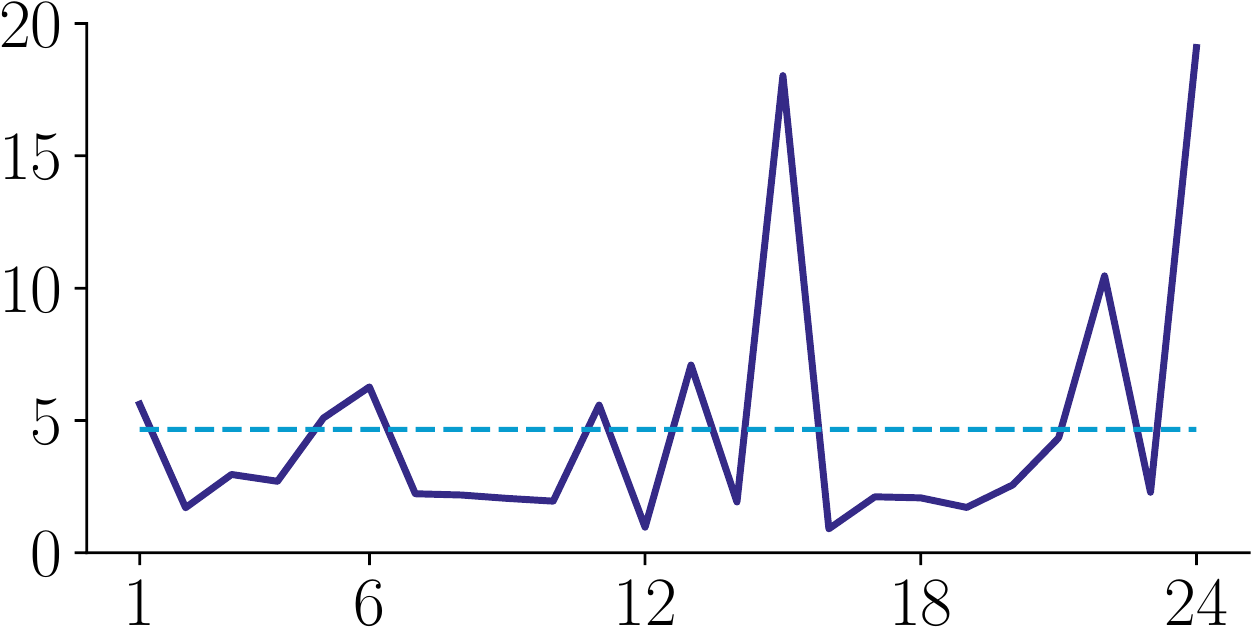}
  \put(10,98){\scriptsize{\textcolor{par_1}{$\boldsymbol{p=3}, $ 24 curves}}}
  \put(10,46){\scriptsize{\textcolor{par_1}{$\boldsymbol{p=3}, $ 24 curves}}}
  \put(-11,5){\scriptsize{\rotatebox{90}{\textbf{CPU time (s)}}}}
\put(40,-7){\scriptsize{\textbf{Curves}}}
\put(36,17){
  \tikz{  \draw[draw=white,fill=white] (0,0) rectangle (2,0.4) }}
  \put(40,19){\scriptsize{\textcolor{par_4}{\textbf{average = 5}}}}
  \end{overpic}
 \end{tabular}
 \end{tabular}
 \caption{\label{fig:nacelle_weak_perf}parallel performance tests approximating the curves of the nacelle model with degrees $p=2,\ 3,$ respectively. Left: weak scaling test using the same curve for every core (1 curve) vs. each core optimizes a different curve (24 curves). Right: iterations and CPU times including average data for the case $p=3,$ 24 curves.}
\end{figure}
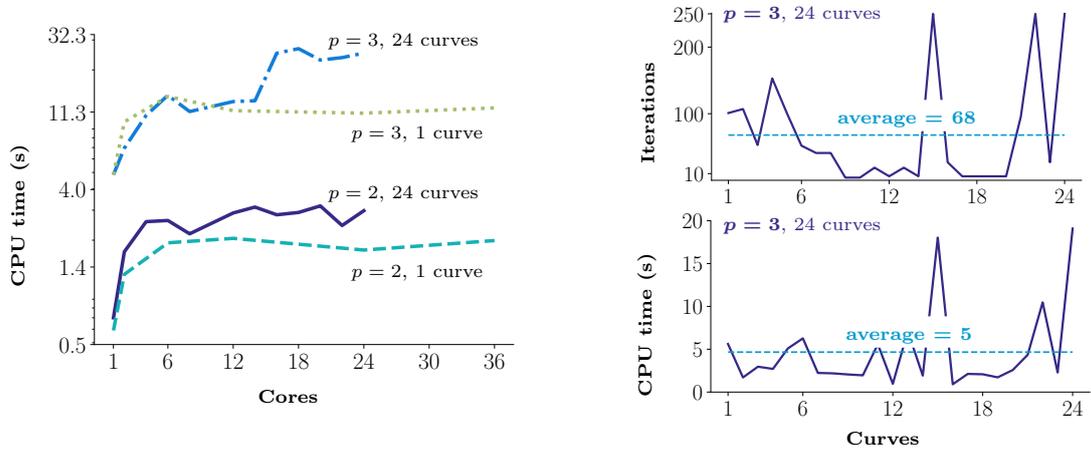

We now test the performance by distributing the curves among the available cores. In practice, since we need to approximate all the curves in the CAD, a strong scaling performance is a better indicator. Moreover, as for the previous study (Figure \ref{fig:element_iters}) where some elements converged with noticeably more iterations, different curve shapes also challenge the optimizer. In Figure \ref{fig:nacelle_weak_perf} (left) we show plots of the CPU times for the nacelle model from Figure \ref{fig:nacelle_fix_free} using the same curve in every core (curve 1) and each core optimizing the different curves (24 curves). This would correspond to a weak scaling study. There is a jump in the CPU times at core 17 for the $p=3$ case. Looking at the right plots, we see that three cores converged in 250 iterations, taking up to 19 seconds of computing time compared to the average of 5 seconds. 

\begin{figure}[h!]
 \centering
 \hspace{10pt}
 \begin{tabular}{cc}
 \begin{overpic}[width = keepaspectratio, width = 0.4\textwidth]{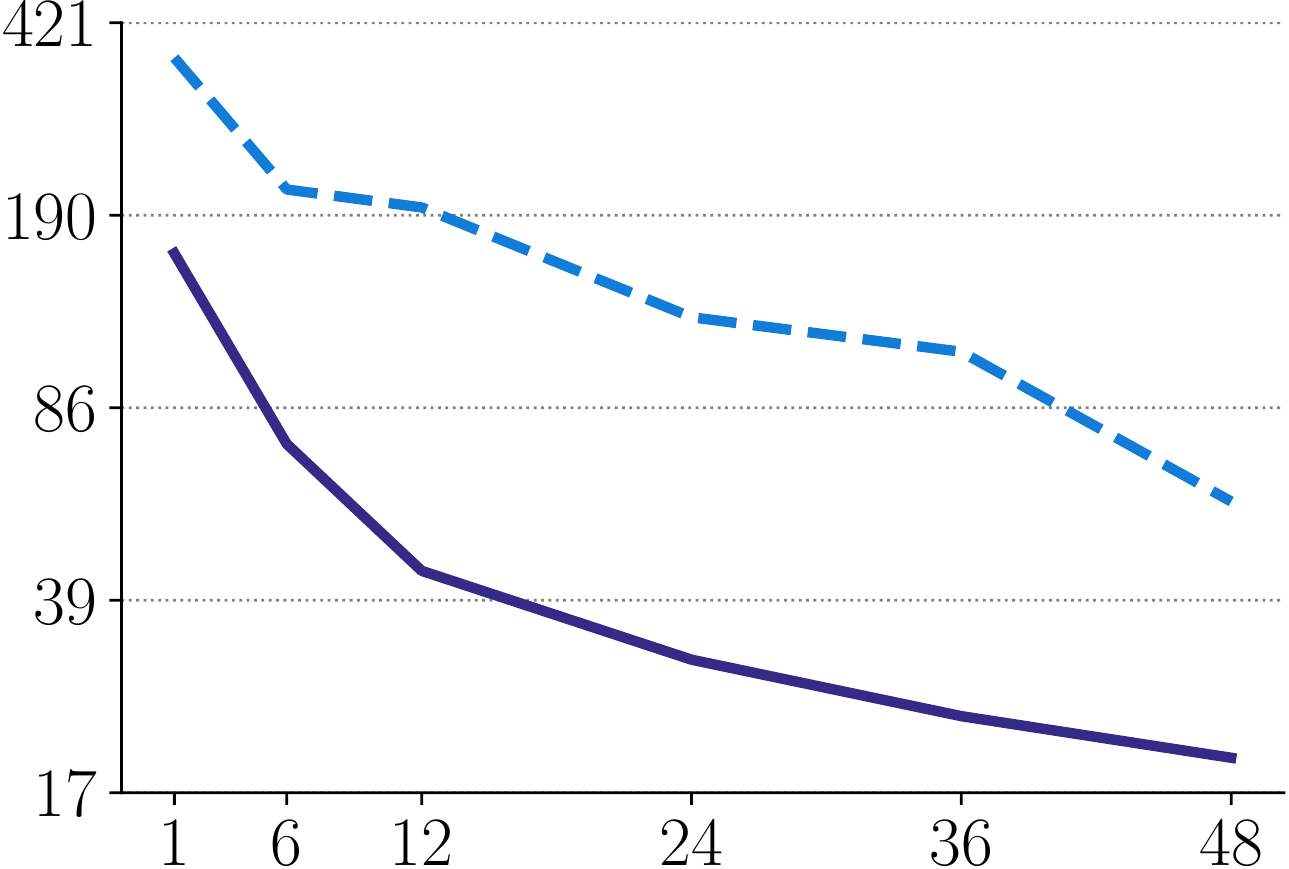}
  \put(-8,18){\scriptsize{\rotatebox{90}{\textbf{CPU time (s)}}}}
  \put(45,-7){\scriptsize{\textbf{Cores}}}
  \put(40,22){\scriptsize{\textcolor{par_1}{  \textbf{Aircraft 1, }$\boldsymbol{p=2}$} }}
  \put(40,52){\scriptsize{\textcolor{par_3}{  \textbf{Aircraft 2, }$\boldsymbol{p=3}$} }}
  \end{overpic}
 & \hspace{30pt}
 \begin{overpic}[width = keepaspectratio, width = 0.4\textwidth]{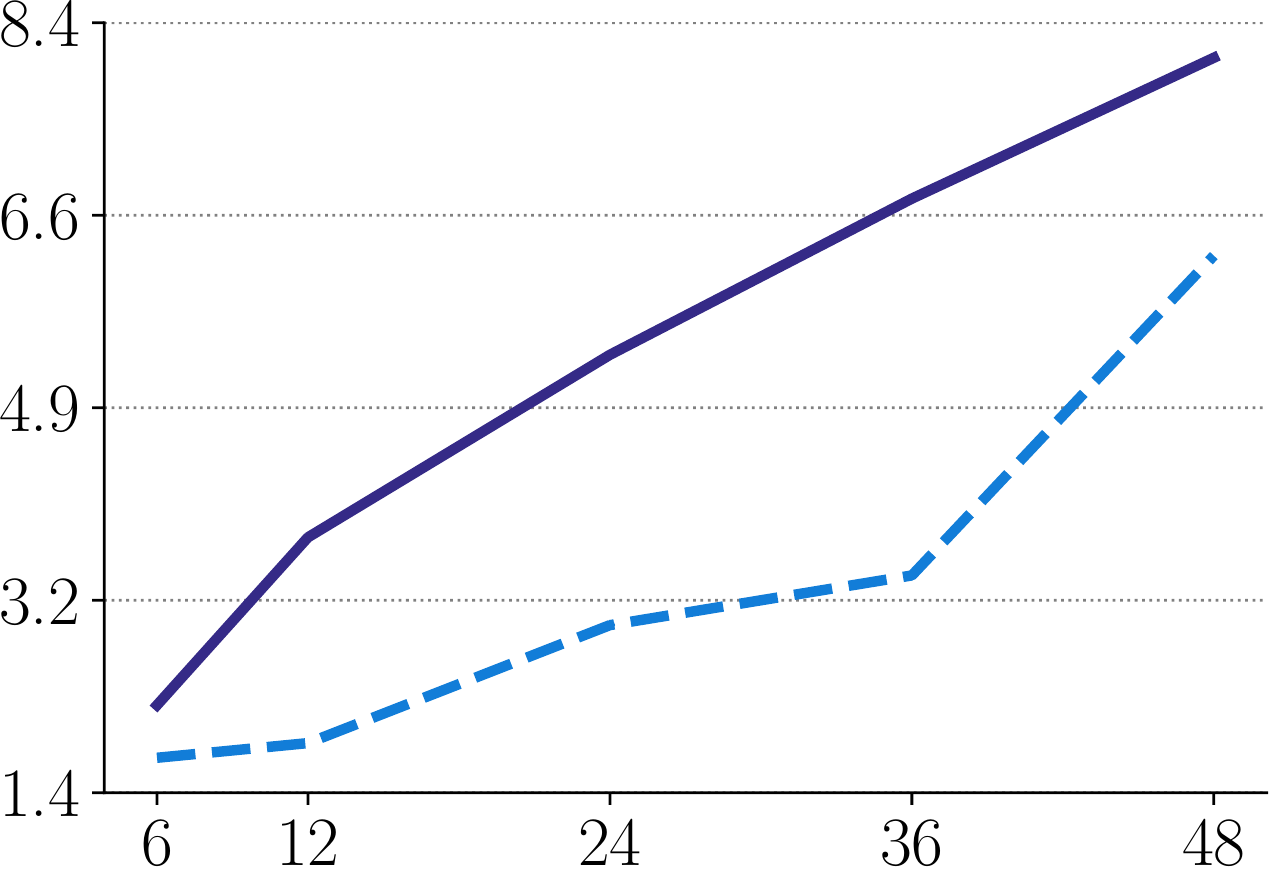}
  \put(-8,18){\scriptsize{\rotatebox{90}{\textbf{Speedup $t_s/t_p$} }}}
  \put(45,-7){\scriptsize{\textbf{Cores}}}
  \put(30,48){\scriptsize{\textcolor{par_1}{ \textbf{A. 1, }$\boldsymbol{p=2}$} }}
  \put(30,22){\scriptsize{\textcolor{par_3}{ \textbf{A. 2, }$\boldsymbol{p=3}$} }}
  \put(90,65){\scriptsize{\textcolor{par_1}{ \textbf{8.1}}}}
  \put(90,50){\scriptsize{\textcolor{par_3}{ \textbf{6.3}}}}
  \end{overpic} 
 \end{tabular}
  \caption{\label{fig:planes_strong} CPU times (left) and speedup (right) increasing the number of cores when optimizing the aircraft models using meshes of degrees $p=2,\ 3,$ respectively.}
\end{figure}
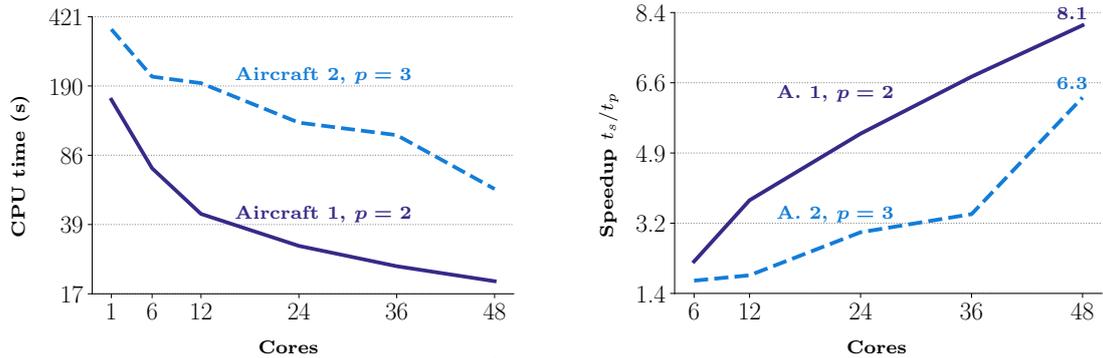

Finally, in Figure \ref{fig:planes_strong} we show the CPU times and speedup factors as we increase the number of cores when optimizing the two aircraft models. Recall that this implies 141 curves for the first model and 86 for the second (see Figure \ref{fig:plane_models}). We can see that the computing times for the $p=2$ case have decreased from 2.7 minutes to 18 seconds and from 6.1 minutes to 58 seconds for the $p=3$ case. 

\section{Conclusions}\label{sec:conclusion}
We have developed a parallel and robust solver that minimizes the constrained disparity measure. The optimization scheme combines a globalized Newton method with the Zhang-Hager nonmonotone line search, and a log barrier penalty term to avoid curve tangling. The constrained disparity is sub-optimal in terms of the error but still yields super-convergence. We have numerically shown how both disparities (original and constrained) are $2p$ super-convergent for 2D curves and  $\lfloor\frac 32(p-1)\rfloor + 2$ for curves in 3D space. The experiments on aircraft prototypes show how the optimized meshes improve the initial approximation error by at least 85\%. 
 
The \emph{Julia} EGADSlite interface combines the use of geometry in parallel with the powerful and relatively simple \emph{Julia} tools for HPC, making this methodology practical for applications in computer clusters. Computing each element in a distributed setting resulted in a low speedup due to few elements limiting the overall performance. On the other hand, optimizing the curves in parallel  resulted in a 8.1 speedup when running on 48 cores. Our results show that we can generate optimal meshes at reasonable computational times. Using 48 cores, it takes 18 seconds to optimize 141 curves of degree $p=2$ and 58 seconds to optimize 86 curves of degree $p=3$.

\section{Acknowledgements}
 This project is funded by the  European  Union's  Horizon  2020  research and innovation programme under the Marie Sk\l{}odowska-Curie grant agreement No 893378. The HPC experiments were run on the \emph{BSC MareNostrum 4}. All the CAD models were borrowed from the \emph{ESP} database. The author would like to thank X. Roca for his encouragement in writing this manuscript,  M. Galbraith for his enormous help with the EGADS \emph{Julia} wrap, and of course, Bob Haimes, for everything.

\bibliography{demoref}

\end{document}